%% file: main.tex
\documentclass[11pt, oneside]{article}

\pdfoutput=1

\usepackage{algorithm, algorithmic, amsmath, amsthm, amssymb, amsbsy, amsfonts, bbm, bm, xcolor, float, footmisc, geometry, graphics, graphicx, hyperref, import, mathrsfs, overpic, stmaryrd, subcaption, svg, upgreek, verbatim, wasysym}

\definecolor{brickred}{rgb}{0.8, 0.0, 0.0}

\DeclareFontShape{OMX}{cmex}{m}{b}{<-> cmexb10}{}
\SetSymbolFont{largesymbols}{bold}{OMX}{cmex}{m}{b}

\newcommand{\E}{\mathbb{E}}
\newcommand{\R}{\mathbb{R}}

\newcommand{\Z}{\mathbb{Z}}
\newcommand{\N}{\mathbb{N}}

\renewcommand{\P}{\mathbb{P}}

\newcommand{\bv}{\mathbf v}

\bmdefine{\bbeta}{\beta}
\renewcommand{\AA}{{\normalfont \textbf{A}}}
\renewcommand{\aa}{{\normalfont \textbf{a}}}
\newcommand{\BB}{{\normalfont \textbf{B}}}
\newcommand{\bb}{{\normalfont \textbf{b}}}
\newcommand{\ee}{{\normalfont \textbf{e}}}

\newcommand{\GG}{{\normalfont \textbf{G}}}
\newcommand{\kk}{{\normalfont \textbf{k}}}

\newcommand{\NN}{{\normalfont \textbf{N}}}

\newcommand{\PP}{{\normalfont \textbf{P}}}
\renewcommand{\tt}{{\normalfont \textbf{t}}}
\newcommand{\uu}{{\normalfont \textbf{u}}}
\newcommand{\vv}{{\normalfont \textbf{v}}}
\newcommand{\ww}{{\normalfont \textbf{w}}}
\newcommand{\xx}{{\normalfont \textbf{x}}}
\newcommand{\XX}{{\normalfont \textbf{X}}}

\newcommand{\Pn}{{\normalfont \text{P}}_n}

\newtheorem{thm}{Theorem}

\newtheorem{lem}{Lemma}
\newtheorem{cor}{Corollary}

\newtheorem{prop}{Proposition}

\title{Fluctuations of lattice zonotopes and polygons}

\author{Théophile Buffière\footnote{LAGA Université Sorbonne Paris Nord.} \thanks{LIPN Université Sorbonne Paris Nord.} , Philippe Marchal $^*$}

\date{}

\begin{document}

\maketitle

\import{Parties/Abstract-Intro/}{abstract.tex}

\section{INTRODUCTION}\label{section_intro}
\import{Parties/Abstract-Intro/}{introduction.tex}

\section{GEOMETRIC AND PROBABILISTIC MODEL}\label{section_model}

\import{Parties/Development/}{1model.tex}
\leavevmode\\

\section{LIMIT THEOREMS IN DIMENSION $d$}\label{section_limit_thm}

\import{Parties/Development/}{2variation.tex}

\leavevmode\\



\section{A DONSKER-LIKE THEOREM IN DIMENSION 2}\label{section_donsker}
\import{Parties/Development/}{3donsker.tex}
\leavevmode\\

\section{BROWNIAN FLUCTUATIONS OF LARGE POLYGONS}\label{section_comment}
\import{Parties/Development/}{5polygons.tex}
\leavevmode\\
\leavevmode\\


\noindent{\Large \textbf{Acknowledgements}} \leavevmode\\ 

\import{Parties/Development/}{8acknowledgment.tex}
\leavevmode\\

\bibliographystyle{amsplain}
\bibliography{bibli}
\end{document}

%% file: Parties/Abstract-Intro/abstract.tex
\begin{abstract}
    Following Barany et al. \cite{Barany-Bureaux:zonotopes}, who proved that large random lattice zonotopes converge to a deterministic shape in any dimension after rescaling, we establish a central limit theorem for finite-dimensional marginals of the boundary of the zonotope. In dimension 2,  for large random convex lattice polygons contained in a square, we prove a Donsker-type theorem for the boundary fluctuations, which involves a two-dimensional Brownian bridge and a drift term that we identify as a random cubic curve.
\end{abstract}

%% file: Parties/Abstract-Intro/introduction.tex
The study of large convex lattice polytopes, starting with the famous
question of their enumeration raised by Arnold \cite{Arnold:polygons}, has been a long-standing
question where, in dimension $\geq 3$, all natural problems regarding
their asymptotic shape remains essentially open. 

In this context, the study of lattice zonotopes, a subclass of convex lattice
polytopes whose definition is given below, seems more tractable.
In particular, an important result by Barany, Bureaux, and Lund \cite{Barany-Bureaux:zonotopes} is that in any dimension, in a given cone,
the shape of large random zonotopes converges to a deterministic limit after rescaling.

Following this result, it is natural to investigate the second order
asymptotics. A convenient way to proceed is to study the fluctuations 
of tangent points of the boundary away from their expected position. We show that this leads to Gaussian limits with a renormalizing factor $n^{\frac{d+2}{2(d+1)}}$.




To be more precise, recall that a zonotope in $\mathbb{R}^d$ is defined as the
Minkowski sum of a finite set $E$ of vectors $\in \mathbb{R}^d$, the vectors in $E$
being called the generators of the zonotope. A zonotope is a lattice zonotope if the vectors have integer coordinates.
 If we restrict ourselves to vectors with nonnegative
 integer coordinates and such that $\sum_{v\in E}v=(n,n\ldots n)$, we get a finite number of
 zonotopes and if we pick one uniformly at random, we have:

\begin{thm}\label{thm:1_CLT_under_Qn}
  Let $Z$ be a random, uniform lattice zonotope starting at the origin, ending at $(n,n\ldots n)$, and with generators in $\N^d$. Let $\uu\in\R^d$ and let $\XX^n_{\uu}$ be the point of the boundary of $Z$  tangent to the hyperplane with normal vector $\uu$, as defined by \eqref{def:X_uu}.
  Then there exists a symmetric matrix $\Gamma_\uu$, given by \eqref{eq:limite_moyenne_covariance} , such that
\begin{equation}\label{CLT1} (n^{-\frac{d+2}{d+1}}\Gamma_{\uu}) ^{-1/2}\left(\XX^n_{\uu} -\E(\XX^n_{\uu})\right) \overset{(d)}{\underset{n\rightarrow \infty}{\longrightarrow}} \mathcal{N}\end{equation}
  where $\mathcal{N}$ is a standard, $d$-dimensional Gaussian variable.
\end{thm}

We point out that, while the boundary of a $d$-dimensional zonotope is a
$d-1$-dimensional object, the fluctuations of the tangent point
$\XX^n_{\uu}$ are $d$-dimensional. 

 If  $H$ is a hyperplane and does not contain any generator of the zonogon, then there are two points of the boundary tangent to $H$, and they are symmetric with respect to the center $(n/2, \ldots n/2)$ of the zonogon. Consequently, the fluctuations of these two points away from their mean are the same. On the other hand, if $H$ does contain some generators of the zonogon, then the set of points of the boundary tangent to  $H$ is the union of two faces of the zonogon, the dimension of these faces being the dimension of the vector space spanned by these generators contained in $H$. Again, these two faces are symmetric with respect to the center of the zonogon. In that case, we have to choose a point on these faces, and we do so using Formula \eqref{def:X_uu}.
 
Let us mention that we have an explicit expression, not only for $\Gamma_{\uu}$ but also for
$\E(\XX^n_{\uu})$. In fact, the result can be generalized to finite-dimensional marginals: if one takes a family
$u_1, \ldots u_k$ of vectors and looks at the tangent points $\XX^n_{u_1}, \ldots \XX_{u_k}$,
then the $k$-tuple $(\XX^n_{u_1}-\E(\XX^n_{u_1}), \ldots \XX^n_{u_k}-\E(\XX^n_{u_k}))$
converges in law to a Gaussian random vector with an explicit covariance structure, see Proposition 4.
Moreover, the result of convergence in distribution can be refined into a local limit theorem. 
See Section \ref{section_limit_thm}.

In the 2-dimensional case, zonogons are just centrally symmetric
convex lattice polygons. Alternatively, any convex
lattice polygon can be viewed as the union of four arcs of zonogons. Moreover,  if we pick a random convex lattice polygon contained in a large square, each of these arcs converges
into an arc of parabola,
as shown in the seminal papers of Barany \cite{Barany:limitshape}, Sinai \cite{Sinai:polygonallines}, and Vershik \cite{Vershik:limit}. In this setting, our result on
fluctuations can be extended to a functional limit theorem as in Donsker's theorem. 

In order to state a rigorous result, let us introduce the following
notation. Consider a random, uniform, convex lattice polygon $\mathcal{P}_n$ contained in the
square $[-n,n]^2$.
Let $\mathcal{A}_n=(A_n, A'_n)$ be the southern-most segment of $\mathcal{P}_n$ and $S_n$
be the ``south pole'', that is, $S_n=(0,-n)$.
Both $A'_n$ and $A_n$ should be close to $S_n$ and we would like to
quantify this more precisely.
Likewise, we can define $B_n, B'n, C_n,C'_n D_n,D'_n$, which should be close respectively to
$E_n,N_n, W_n$ where $E=(n,0)$ etc. See Figure 1.
Finally, we denote by $X:\R^2\to\R$, resp. $Y$, the projection on the first
(resp. second) coordinate. Then we get the following result:

\begin{figure}[H]
    \centering
    \begin{overpic}[width = 16cm]{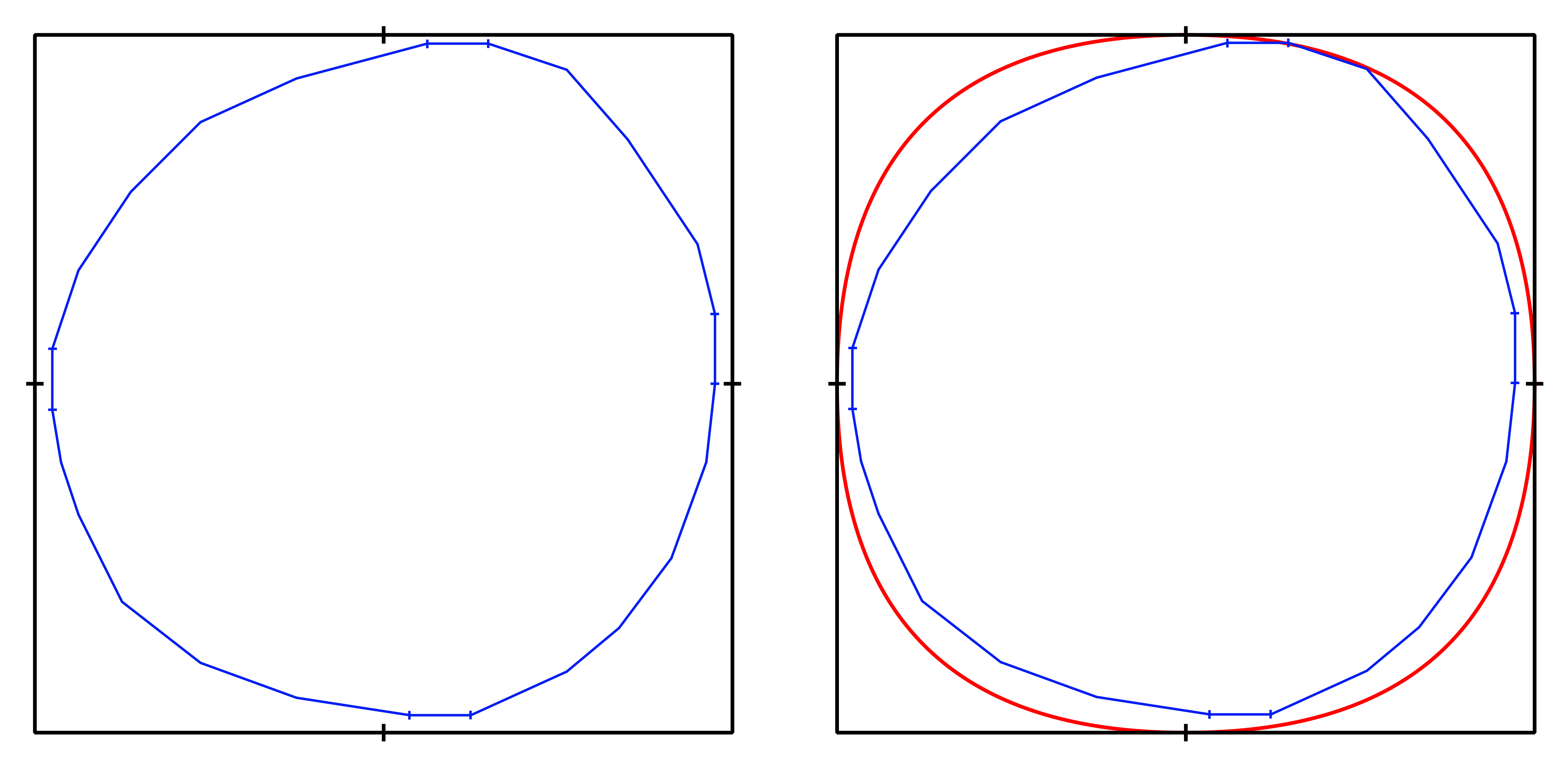}
        \put(-2,23.6){$W_n$}
        \put(4,25.5){\color{blue}$D_n$}
        \put(4,21.4){\color{blue}$D'_n$}
        \put(23,-0.5){$S_n$}
        \put(24,4){\color{blue}$A_n$}
        \put(27.5,4){\color{blue}$A'_n$}
        \put(23,47.5){$N_n$}
        \put(26,43.5){\color{blue}$C_n$}
        \put(29,43.5){\color{blue}$C'_n$}
        \put(47.5,23.6){$E_n$}
        \put(42,28){\color{blue}$B'_n$}
        \put(42,23.6){\color{blue}$B_n$}
    \end{overpic} %
    \caption{A random polygon in the square $[-n,n]^2$, where the extremal segments are labeled by their extremities, and the limit shape shown by Barany \cite{Barany:limitshape}, Sinai \cite{Sinai:polygonallines}, and Vershik \cite{Vershik:limit}.}%
    \label{fig:forme_limit_polygon}%
\end{figure}
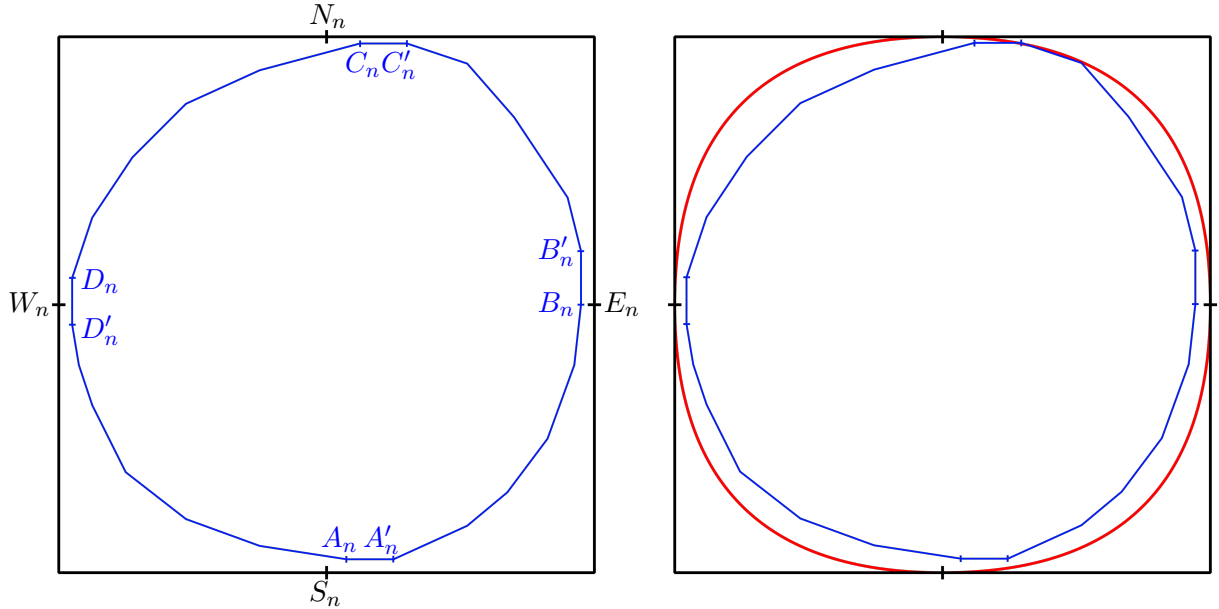


\begin{thm}\label{thm:brownian_polygone}
(i) The quadruple $$n^{-2/3}(Y(A_n- S_n), X(B_n- E_n),Y(C_n- N_n),X(D_n- W'_n))$$
converges in probability to $\delta_{(0,0,0,0)}$.

(ii)
The quadruple $$n^{-2/3}(X(A_n- S_n), Y(B_n- E_n),-X(C_n- N_n)-Y(D_n- W_n))$$
converges in distribution to a Gaussian random variable $(R,S,T,U)$ with
density
$$
C_1\exp\left(-\frac{1}{18}[(r-s)^2+(s-t)^2+(t-u)^2+(u-r)^2]-\frac{1}{3}
[r^2+s^2+t^2+u^2]
\right)
$$
for some normalizing constant $C_1>0$.

(iii) For every $t\in [0,1]$, let $\XX_n(t)$ be the point of the boundary of $\mathcal{P}_n$ with negative $y$-coordinate and tangent to
the vector $(t,1-t)$. Let $\overline{\XX}_n(t)=\E(\XX_n(t))$.
Then there exists a  continuous family of nonsingular matrices $(Q(t), 0\leq t\leq 1)$ such that for all $r,s\in \R$,
denoting for each $n$ the event
$$\mathcal{E}_n(r,s)=\{\lfloor n^{-2/3}(X(A_n- S_n)\rfloor =r,\ \lfloor n^{-2/3}(Y(B_n- E_n)\rfloor =s\}$$ 
we have the convergence
of conditional processes
$$\left(n^{-2/3}(\XX_n(t)-\overline{\XX}_n(t))|\mathcal{E}_n(r,s)), t\in[0,1]\right)
\stackrel{(d)}{\to} \left(\begin{matrix}
                r\\
                0
            \end{matrix}\right)+\left(\boldsymbol{\mu}_{r,s}(t) +Q(t)\beta_t, t\in[0,1]\right)$$
where $(\beta_t)$ is a standard 2-dimensional Brownian bridge and  $\boldsymbol{\mu}_{r,s}$
is a cubic curve parameterized by
$$
         \boldsymbol{\mu}_{r, s}(t) = \left(\begin{matrix}
                -2t(t-1)^2  & t(2t^2 -5 t + 4)\\
                t^2(2t - 1) & -2 t^2(t-1)
            \end{matrix}\right)\left(\begin{matrix}
                s\\
                -r
            \end{matrix}\right)
$$
\end{thm}

Of course, (iii) is only stated for one of the four arcs of the polygon but the result is true for each arc. Note that $\XX_n(0)=A_n$ and $\XX_n(1)=B_n$.
As in Theorem 1, the tangent point in Theorem 2 (iii) is defined by \eqref{def:X_uu}. Anyway, it follows from \cite{Buffiere:number} that the edges of $\mathcal{P}_n$ have length of order $n^{\frac{1}{3}}$. Therefore, even if the set of tangent points is a whole edge, we could choose any point on this edge as $\overline{\XX}_n(t)$ and because of the renormalizing factor $n^{-2/3}$, this would not change the result. For the same reason, we could  replace $(A_n, B_n,C_n,D_n)$ with $(A'_n, B'_n,C'_n,D'_n)$ in the statement of the theorem.

We could re-express (iii) by saying that $(Q(t)\beta_t)$ is a 2-dimensional Gaussian process with a continuous
family of covariance matrices $(Cov_t)$.  It follows from the computations in Section 4 that there exists a continuous family of orthogonal matrices $(\mathcal{O}_t)$, with $\mathcal{O}_0=Id$, such that
\begin{align}
    \mathcal{O}_t Cov_t = \left(\frac{\zeta(2)}{\zeta(3)}\right)^{1/3}  \begin{pmatrix}
2\left(4t^2 - 6t +3 \right)(t-1)^3 t & \left( 8t^2 - 8t +3 \right) (t-1)^2 t^2\\
\left( 8t^2 - 8t +3 \right) (t-1)^2 t^2  & -2 \left((4t^2-2t+1\right)t^3(t-1) 
\end{pmatrix}
\end{align}
In particular, for small  $t$, the fluctuations of the process are of order
$t^{1/2}$ in the $x$-coordinate and  $t^{3/2}$ in the $y$-coordinate. We have similar estimates for $t$ close to 1.

As a consequence, we could write informally that if  $t\asymp n^{-\frac{1}{3}}$,
$$
n\overline{\XX}_n(t)\asymp \left(\begin{matrix}
               n^{\frac{2}{3}} \\
               n^{\frac{1}{3}} 
            \end{matrix}\right)\ ,\ \
n^{\frac{2}{3}}\E|Q(t)\beta_t|\asymp \left(\begin{matrix}
               n^{\frac{1}{2}} \\
            n^{\frac{1}{6}} 
            \end{matrix}\right)\ ,\ \ 
n^{\frac{2}{3}}\mu_{r,s}(t)\asymp \left(\begin{matrix}
              n^{\frac{1}{3}} \\
             1 
            \end{matrix}\right)             
$$

where $a_n \asymp b_n$ means that there exist two positive constants $c<C$ such that
for each $n$, $cb_n<a_n<Cb_n$.

Note that in (iii), we state a conditional result. If we average over the law of $X(A_n- S_n), Y(B_n- E_n)$, we find that the mean of the asymptotic cubic curve is zero. That is, if we choose $(R,S,T,U)$ according to the Gaussian distribution given in (ii), 
then for every $t\in[0,1]$, $\E(\boldsymbol{\mu}_{R, S}(t))=0$.

More details on  $\boldsymbol{\mu}_{r, s}$, in particular the proof that it is a cubic curve, are also given in Section 4, see Proposition 6, where we have to identify 
$$\mu_{r,s}(t)=\begin{pmatrix}
0 & 1\\
1 & 0
\end{pmatrix}\nu_{-r,s}(t)$$ 
The curve $\mu_{r,s}$ has a cusp if $-r/s\notin[1/2,2]$. In particular, suppose that $r<0$, $s>0$ and $-r/s>2$. Then $\boldsymbol{\mu}_{r, s}$ starts at $(0,0)$, ends in the positive quadrant, namely at $(-r,s)$,
and yet for $t<(-r-2s)/3(-r-s)$, both coordinates of the speed $\boldsymbol{\mu}_{r, s}'(t)$
are negative. This may seem counter-intuitive.

The fact that for fixed $r,s$, the curve  $\boldsymbol{\mu}_{r, s} $ is cubic is in sharp
contrast with the usual situation where the drift is linear: if $(B_t, 0\leq t \leq 1)$ is a Brownian motion started at 0 and conditioned to end at $a$, then $(B_t)$ has the form 
$B_t=at+\beta_t$ where $(\beta_t)$ is a Brownian bridge. In dimension 1, Brownian motion with a parabolic drift has been widely studied in connection with various problems such as statistical estimators,  random partitions,  epidemics models, Burgers turbulence etc. See for instance \cite{janson:maxbrownian10,Joseph:randomgraph,vanderhofstad:parabolic,GIRAUD:200341,Groeneboom:BrownianMW89}. On the other hand, we are not aware of other instances of a cubic drift in the literature.

Let us also mention that other models of random polytopes have been studied in the literature. In particular, Calka and coauthors \cite{CALKA:rdpolytope, calka:brownianlimit} consider the convex hull of Poisson point processes whose intensity goes to infinity, thereby obtaining random polytopes that fill nearly all the space available, and look at the fluctuations of the boundary. This, however, is very different from our setting where the geometry of the underlying lattice plays a significant role and therefore, their results and ours are of a different nature.

The remainder of this paper is organized as follows.
We introduce our main tools in the next section. In particular, while the two results stated in the introduction deal with random zonotopes or polygons with uniform distribution, a key tool we will use is Boltzmann distributions, which are defined in this part of the paper. We state and prove our limit theorems in
dimension $d$ in Section \ref{section_limit_thm}. Section \ref{section_donsker} is devoted to the functional, Donsker-like result on
zonogons in dimension 2. Finally, Section \ref{section_comment} extends these results to polygons and gives the proof of Theorem 2.

%% file: Parties/Development/1model.tex
The results of this paper are based on the limit theorems proved in Section \ref{section_limit_thm}, which are dealing with zonotopes in cones of $\R^d$, for any dimension $d$. Thereafter, the cone $\mathcal{C}$ is a closed convex salient and pointed cone, which means that $0\in \mathcal{C}$, no pair $(\xx, -\xx)$ lies in $\mathcal{C}$ for some non zero vector $\xx$, and its interior $\text{int } \mathcal{C}$ is not empty. Let $\kk$ be vector of $\Z^d \cap \mathcal{C}$. The geometric and probabilistic models are the same as in \cite{Barany-Bureaux:zonotopes}; we tried to keep the same notation when possible to ease the readability of our results. 

\subsection{Zonotopes in cones}

A \emph{zonotope} is a convex geometric object defined as the Minkowski sum of $k$ segments, called its \emph{generators}. More specifically, an \emph{integral zonotope} $Z$ is a polytope for which there exist $k \in \N$ and $\vv_1, ..., \vv_k \in \Z^d$ such that

\begin{align*}
    Z &= \left\{\sum_{i=1}^k \alpha_i \vv_i  \:\:| \:\: (\alpha_1,...,\alpha_k) \in [0,1]^k\right\}\\
    &= \text{Conv}\left\{ \sum_{i=1}^k \epsilon_i \vv_i  \:\:| \:\: (\epsilon_1,...,\epsilon_k) \in \{0,1\}^k  \right\}
\end{align*}

\begin{figure}[H]
    \centering
    \subfloat[\centering 2 dimensions ] {{ \includegraphics[width=6cm]{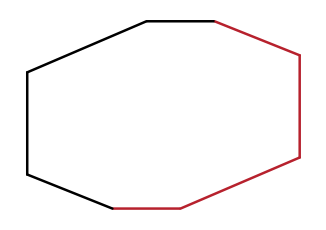} }}%
    \qquad
    \subfloat[\centering 3 dimensions] {{ \includegraphics[width=6cm]{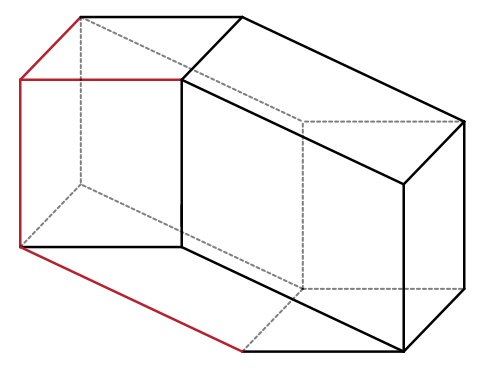} }}%
    \caption{Zonotopes made from 4 red-colored generators.}%
    \label{fig:example_zono}%
\end{figure}

The endpoint of $Z$ is $\sum_{i = 1}^k \vv_i$. We define $\mathcal{Z}(\mathcal{C}, \kk)$ as the set of integral zonotopes in $\mathcal{C}$ that end at $\kk$. This set is clearly finite as the generators of the zonotopes of $\mathcal{Z}(\mathcal{C}, \kk)$ are in $\Z^d$.\medskip

Given a multiset (an unordered set of elements with repetition allowed) $V = \{\vv_1,..,\vv_k \}\in (\Z^d\setminus\{\boldsymbol{0}\})^k$ and the integral zonotope $Z$ determined by $V$, the set of generators uniquely defines a zonotope but the converse is not true.  To determine a multiset uniquely from an integral zonotope, we define the set $\P_d$ of primitive vectors of $\Z^d$. A vector $(z_1,...z_d) \in \Z^d$ is primitive if $\gcd(z_1, ..., z_d) = 1$, and therefore notice that $0\notin \P_d$ (we prefer the notation $\P_d$ rather than $\P^d$ because $\P_d$ is very different from $ (\P_1)^d$). There is a unique multiset $W = \{ \ww_1, ..., \ww_l\} \in \P_d^l$ of elements of $\P_d$ that determines $Z$, constructed that way: given a generator $\vv \in \Z^d$ of $Z$, $\vv$ can be uniquely written as $m_i \ww$, with $\ww \in  \P_d$. Then add $m_i$ copies of $\ww$ in $W$. For two different multisets $U$ and $V$ that determine $Z$, $U$ and $V$ give the same $W$ with the construction above. \medskip

We have now a one-to-one correspondence between integral zonotopes in $\mathcal{C}$ and finite multisets of $\P_d \cap \mathcal{C}$. We denote the \textit{endpoint} $\kk$ of $Z$ (and we say that $Z$ ends at $\kk$) as 
$$ \kk = \sum_{i =1}^{l} \ww_i.$$

A multiset of $\P_d \cap \mathcal{C}$ describing a zonotope such as $W$ is called a \textit{strict integer partition} (see Section 2 in \cite{Barany-Bureaux:zonotopes}) of the vector $\kk$ from $\mathcal{C}$. Notice that for a given $\kk \in  \Z^d \cap \mathcal{C}$, any strict integer partition of $\kk$ from $\mathcal{C}$ is finite and can be encoded as a function with finite support. Indeed, for a strict integer partition $W$, we define the function of multiplicities $\omega : \P_d \cap \mathcal{C} \rightarrow \Z_+ $ where $\omega(\xx) = \text{card}\{ i \in \{1,...,l \}\: |\: \ww_i = \xx \}$ for any $\xx \in  \P_d \cap \mathcal{C}$. \medskip

Let $\Omega(\mathcal{C})$ be the set of nonnegative integer-valued functions $\omega : \P_d \cap \mathcal{C} \rightarrow \Z_+$ with finite support. There is a one-to-one correspondence between $\Omega(\mathcal{C})$ and the set of integral zonotopes in $\mathcal{C}$. Additionally, we naturally define the endpoint of $\omega$ 
\begin{align}\label{eq:X(omega)}
    \XX(\omega) = \sum_{\xx \in \P_d \cap \mathcal{C}} \omega(\xx) \xx.    
\end{align}

Picking a random zonotope in $\mathcal{Z}(\mathcal{C}, \kk)$ is equivalent to picking $\omega \in \Omega(\mathcal{C})$ such that $\XX(\omega) = \kk$. In the sequel, $Z(\omega)$ denotes the zonotope that corresponds to $\omega$, and $\omega(Z)$ the element of $\Omega(\mathcal{C})$ that corresponds to $Z$. 


\subsection{The probabilistic model}

Fix $\kk \in \text{int }\mathcal{C} \cap \Z^d$. We will use the correspondence between integral zonotopes and strict integer partitions described above to define the probability distribution on the set of zonotopes. For $\omega \in \Omega(\mathcal{C})$, we denote $Z(\omega)$ the associated zonotope. \medskip    

We define the probability distribution $\Pn$ for all $\omega \in \Omega(\mathcal{C})$, depending on the parameters $\beta_n \in \R_+$ and $\aa \in \Z^d \cap \text{int}\mathcal{C}$ (respectively made explicit in (\ref{eq:beta_n}) and in (\ref{eq:def_a})) by 

$$\Pn (\omega)= \frac{1}{Z_n(\aa)} e^{-\beta_n \aa \cdot \XX(\omega)}, \:\:\:\:\: \text{ where } Z_n(\aa) = \sum_{\omega \in \Omega(\mathcal{C})} e^{-\beta_n \aa \cdot \XX(\omega)}. $$ 

This probability distribution is known as the Boltzmann probability distribution, as it is directly inspired by the Boltzmann distribution in statistical physics. $Z_n(\aa)$ is called the \textit{partition function} of the model. In the sequel, we fix $\beta_n$ throughout the paper to be

\begin{align}\label{eq:beta_n}
    \beta_n = \sqrt[d+1]{ \frac{\zeta(d+1)}{\zeta(d) n }}
\end{align}

The key point is that $\Pn(\omega)$ only relies on the endpoint $\XX(\omega)$, hence two zonotopes ending at the same point have the same probability, and in particular, we define the uniform distribution $\text{Q}_n$ on the elements $\Omega(\mathcal{C})$ ending at $n \kk$ by
$$ \text{Q}_{n \kk} (\omega) = \Pn \left( \omega \: | \: \XX(\omega) = n \kk \right)$$

The parameters $\aa$ and $\beta_n$ are determined in order for $\Pn$ to be close to $\text{Q}_{n \kk}$ when $n$ grows large. The point of using $\Pn$ to approximate $\text{Q}_{n\kk}$ is that $\Pn$ has a much simpler structure. Remind the definition of $\XX(\omega)$ in \ref{eq:X(omega)} as a sum over $\P_d \cap \mathcal{C}$; therefore the exponential becomes

$$ e^{-\beta_n \aa \cdot \XX(\omega)} = \prod_{\xx \in \P_d \cap \mathcal{C}} e^{- \beta_n \aa \cdot \omega(\xx) \xx}. $$

This product structure is passed along to $Z_n(\aa)$ and $\Pn$:

\begin{align}\label{eq:P_n_definition}
    Z_n(\aa) = \prod_{\xx \in \P_d \cap \mathcal{C}} \frac{1}{1 - e^{-\beta_n \aa \cdot \xx}}, \:\:\:\:\: \text{ and } \Pn( \omega) = \prod_{\xx \in \P_d \cap \mathcal{C}} e^{- \beta_n \aa \cdot \omega(\xx) \xx} \left( 1 - e^{- \beta_n \aa \cdot \xx} \right)
\end{align}

We deduce that $(\omega(\xx))_{\xx \in \P_d \cap \mathcal{C}}$ is a mutually independent set under $\Pn$, and that the variable $\omega(\xx)$ has a geometric distribution of parameter $1 - e^{- \beta_n \aa \cdot \xx}$. The simplicity of $\Pn$ lies in the fact that under $\Pn$, a random zonotope is a product of geometric independent variables for each possible primitive generator. \medskip

\subsection{Useful formulas about cones.}

This part is dedicated to recalling a few results from B\'ar\'any Bureau and Lund \cite{Barany-Bureaux:zonotopes} about cones, more precisely Theorem 3.1, and Proposition A.1. We keep the $d$-dimensional cone $\mathcal{C}$ and $\kk \in \Z^d \cap \mathcal{C}$ from above. For a given $\aa \in \Z^d$, we denote the section $\mathcal{C}(\bb =t)$ and the cap $\mathcal{C}(\bb \leq t)$ by 
\begin{align*}
    \mathcal{C}(\bb =t) = \left\{\xx \in \mathcal{C} | \xx \cdot \bb = t \right\}, \\
    \mathcal{C}(\bb \leq t) = \left\{ \xx \in \mathcal{C} | \xx \cdot \bb \leq t \right\}.
\end{align*}

The following proposition (Theorem 3.3 in \cite{Gigena:cones}) fixes $\aa$ in the distribution $\Pn$:

\begin{prop}\label{prop:centre_cone}
    Given $\mathcal{C}$ and a vector $\vv \in \textup{int } \mathcal{C}$, there is a unique $\aa' = \aa'(\mathcal{C}, \vv)$ such that $\aa'$ is the center of gravity of the section $\mathcal{C}(\vv=1)$ and that the cap $\mathcal{C}(\aa')$ is the unique cap that has minimal volume among all caps of $\mathcal{C}$ that contains $\vv$.
\end{prop}

It immediately follows that $\frac{d}{d+1}\aa'$ is the center of gravity of $\mathcal{C}(\vv \leq 1)$, and therefore 
$$ \frac{d}{d+1}\aa' = \frac{1}{\text{Vol} \mathcal{C}(\vv \leq 1)} \int_{\mathcal{C}(\vv \leq 1)} \xx d\xx.$$

In the probability distribution $\Pn$, given $\aa'(\mathcal{C}, \kk)$ obtained by Proposition \ref{prop:centre_cone}, we define $\aa \in \Z^d \cap \mathcal{C}$ the unique vector $\aa = \lambda \aa'(\mathcal{C}, \kk)$ for some $\lambda > 0$ such that
\begin{align}\label{eq:def_a}
    (d+1)! \int_{\mathcal{C}(\aa\leq 1)}\xx d\xx = \kk.
\end{align}

This value will make sense in the proof of Theorem \ref{thm:CLT}. The next proposition is a version of the density of the primitive vectors for cones and homogeneous functions, and will be used hereafter for calculations on $\Pn$:

\begin{prop} \label{prop:equivalent_int_sum}
For $d> 1$, let $f: \mathcal{C} \rightarrow \R$ be a continuously differentiable homogeneous function of degree $h$, i.e. $f(\lambda \xx) = \lambda^h f(\xx)$ for all $\xx \in \mathcal{C}$ and $\lambda \geq 0$. For every $\aa \in \R^d$ such that $\aa \cdot \xx > 0 $ for all $\xx \in \mathcal{C}$ (that is, an element in the interior of the dual cone):

$$ \beta^{d+h} \sum_ {\xx \in \mathcal{C}\cap \mathbb{P}_d} f(\xx) e^{- \beta \aa \cdot \xx} = \frac{1}{\zeta(d)} \int_{\mathcal{C} } f(\xx) e^{- \aa \cdot \xx} d\xx + O(\beta)$$

\end{prop}

This is Proposition A.1 in \cite{Barany-Bureaux:zonotopes}. It is written under the assumption $d\leq 3$, but their proof still stands for $d=2$. At some point hereafter, we will need to highlight that the term $O(\beta)$ is independent of $\mathcal{C}$. For this purpose, we highlight some elements of the proof of Proposition \ref{prop:equivalent_int_sum} that give the following corollary:

\begin{cor}\label{cor:homogeneous_bound_on_cone}
For $d  \geq 2$, if $f$ is a continuously differentiable homogeneous function of degree $h$ and $\mathcal{C} \subset \R_+^d$:

$$ \beta^{2+h} \sum_ {\xx \in \mathcal{C}\cap \mathbb{P}_2} f(\xx) e^{- \beta \aa \cdot \xx} \leq \frac{1}{\zeta(d)} \int_{\mathcal{C} } f(\xx) e^{- \beta \aa \cdot \xx} d\xx + \epsilon(\aa)\beta$$

where $\epsilon(a)$ is independent of $\mathcal{C}$.
\end{cor}

\begin{proof}
Recall that $a \in \text{Int }\R_{+}^d$, as it is proportional to the center of gravity of $\mathcal{C}(\uu \leq 1)$. 
Let $A$ be a compact subset of $\text{Int }\R_{+}^d$ (which is the dual of $\R_{+}^2$, $\{ \vv \in \R^d, \forall \xx \in \R_{+}^d, \vv \cdot \xx > 0 \}$).
Since $A$ is compact and since $\mathcal{C}(\uu \leq 1) \subset \{\xx \in \R_{+}^d, \uu \cdot \xx > 0 \} $, there exists $L>0$ such that for any cone $\mathcal{C}$, $\mathcal{C}(\uu \leq 1)$ is contained in $\left[0, L\right]^d$. Therefore the same arguments as in the proof of  Corollary A.3 in \cite{Barany-Bureaux:zonotopes} lead to the existence of a constant $\gamma_{L, f} >0$ such that 

$$ \underset{\aa \in A}{\sup} \left| \sum_ {\xx \in \mathcal{C}\cap \mathbb{P}_2} f(\xx) e^{- \beta \aa \cdot \xx} - \frac{1}{\zeta(d)} \int_{\mathcal{C} } f(\xx) e^{- \beta \aa \cdot \xx} d\xx \right| \leq \frac{\gamma_{L,f}}{\beta^{2+h}}.$$

This constant $\gamma_{L,f}$ is independent of $\mathcal{C}$. 
\end{proof}

The following relation is widely used in this paper to compute integral values.

$$\int_{\mathcal{C}} e^{- \aa \cdot \xx} d\xx = (d+1)! \int_{\mathcal{C}( \aa \leq 1) } \xx d\xx ,$$

%% file: Parties/Development/2variation.tex
This section improves the results of \cite[Section 4.2]{Barany-Bureaux:zonotopes} on the asymptotic behavior of random zonotopes after rescaling. First, we extend the central limit theorem to any point of the boundary of the zonotope tangent to a given hyperplane. Then we prove a local limit theorem for these points. 
Let $\uu$ and $\vv$ be two vectors of $\mathbb{R}^d$, and  $Z$ a zonotope in $\mathcal{Z}(\mathcal{C},n\kk)$ drawn under $\Pn$ or $\text{Q}_{n\kk}$. Denote  $H_{\uu}$, respectively $H_{\vv}$,  the hyperplane of normal vector $\uu$, respectively $\vv$. We define $\XX^n_{\uu}$ the furthest point from the origin of the face of $Z$ tangent to $H_{\uu}$ that is maximizing the scalar product with $\uu$. Namely, 

\begin{align}\label{def:X_uu}
    \XX^n_{\uu} = \underset{}{\max}\left(\underset{\xx \in Z}{\text{arg max}} (\xx \cdot \uu)\right).
\end{align}

We similarly define $\XX^n_{\vv}$. We can reformulate the definition of $\XX^n_{\uu}$ based on the structure of a zonotope: $\XX^n_{\vv}$ is the endpoint of the zonotope generated by the generators $\xx$ of $Z$ such that $\xx \cdot \uu \geq 0$.\medskip

\subsection{Central Limit Theorem for $\XX_{\uu}$}


 Due to the independence between the generators under $\Pn$, the point $\XX^n_{\uu}(\omega)$ is distributed according to:

$$ \text{P}_{n,\uu}(\omega)  = \prod_{\xx \in \P_d \cap \mathcal{C}(\uu \geq 0)} e^{-\beta_n \aa \cdot \omega(\xx) \xx } (1 - e^{- \beta_n \aa \cdot \xx}) = \frac{1}{Z_{n,\uu}(\aa)} e^{- \beta_n \aa \cdot \XX^n_{\uu}(\omega)},$$

where $Z_{n,\uu}(\aa)$ is the partition function of $\text{P}_{n,\uu}$. For exponential distributions like $\text{P}_{n, \uu}$ and $\Pn$, the expectation and the covariance matrix are known to be written in terms of the derivative of the logarithm of this function, that is for the expectation $\mu_\uu ^n$ and the covariance $\Gamma_\uu ^n$ of $\XX_{\uu}^n$ : 
\begin{align*}
    \mu_\uu ^n = \mathbb{E}_n\left[\XX_\uu ^n \right] = - \nabla \log Z_{n,\uu} (\aa) \:\:\: \text{ and }\:\:\: \Gamma_\uu ^n = \text{Cov} (\XX_\uu ^n) =  \nabla^2 \log  Z_{n,\uu} (\aa).  
\end{align*}

Actually, $\mathcal{C}(\uu\geq 0)$ is a $d$-dimensional cone, and the theorems proved are still true for more complex cones. More generally, given a $d$-dimensional closed convex salient and pointed cone $\mathcal{C}_1$ such that $\mathcal{C}_1 \subset \mathcal{C}$, we define $\XX^n(\mathcal{C}_1)$, $\text{P}_{n, \mathcal{C}_1}$, $\boldsymbol{\mu}^n(\mathcal{C}_1)$, and $\Gamma^n(\mathcal{C}_1)$ analogously:

\begin{align}\label{def:partial_Pn}
    \XX^n(\mathcal{C}_1) = \sum_{\xx \in \P_d \cap \mathcal{C}_1} \omega(\xx) \xx \:\:\: &\text{ and } \:\:\: \text{P}_{n, \mathcal{C}_1} =  \prod_{\xx \in \P_d \cap \mathcal{C}_1} e^{-\beta_n \aa \cdot \omega(\xx) \xx } (1 - e^{- \beta_n \aa \cdot \xx}) \nonumber\\
    \boldsymbol{\mu}^n(\mathcal{C}_1) = - \nabla \log \left( \prod _{\xx \in \P_d \cap \mathcal{C}_1} \frac{1}{1 - e^{ - \beta_n \aa \cdot \xx }} \right) \:\:\: &\text{ and } \:\:\: \Gamma^n(\mathcal{C}_1) = \nabla^2 \log \left(\prod _{\xx \in \P_d \cap \mathcal{C}_1} \frac{1}{1 - e^{ - \beta_n \aa \cdot \xx }} \right)
\end{align}

In particular, we have $\boldsymbol{\mu}^n(\mathcal{C}(\uu \geq 0) ) = \boldsymbol{\mu}^n_{\uu}$ and  $\Gamma^n(\mathcal{C}(\uu \geq 0 )) = \Gamma^n_{\uu}$.  Denote also the rescaled limits $\boldsymbol{\mu}(\mathcal{C}_1)$ (with   $\boldsymbol{\mu}(\mathcal{C}(\uu \geq 0)) = \boldsymbol{\mu}_\uu$ ) and $\Gamma(\mathcal{C}_1)$ (with  $ \Gamma(\mathcal{C}(\uu \geq 0))  =\Gamma_\uu$) respectively as:

\begin{align}\label{eq:limite_moyenne_covariance}
    \boldsymbol{\mu}(\mathcal{C}_1) = (d +1)! \int_{\mathcal{C}(\aa \leq 1) \cap \mathcal{C}_1} \xx d\xx, \:\:\:\:\: \Gamma(\mathcal{C}_1) = \left(\frac{\zeta(d)}{\zeta(d+1)}\right)^{\frac{1}{d+1}} (d+1)! \int_{\mathcal{C}( \aa \leq 1) \cap \mathcal{C}_1} \xx \xx^\intercal d\xx
\end{align}

\begin{prop}[Central limit theorem]\label{thm:CLT}
Let $Z \in \mathcal{Z}(\mathbb{R}^2, n\boldsymbol{1})$ be a random zonotope drawn under $\Pn$,  $\mathcal{C}_1 \subset \mathcal{C}$ be a $d$-dimensional cone with $0 \in \mathcal{C}_1$, and $\XX^n(\mathcal{C}_1)$ be the endpoint of the generators of $Z$ in $\mathcal{C}_1$. Then $\XX^n (\mathcal{C}_1)$ satisfies a central limit theorem in the sense that
$$ (\Gamma^n(\mathcal{C}_1)) ^{-1/2}\left(\XX^n(\mathcal{C}_1) -\boldsymbol{\mu}^n(\mathcal{C}_1) \right) \overset{(d)}{\underset{n\rightarrow \infty}{\longrightarrow}} \mathcal{N}(\boldsymbol{0}, I_d),$$
with $ \underset{n\rightarrow \infty}{\lim} \frac{1}{n} \boldsymbol{\mu}^n(\mathcal{C}_1) = \boldsymbol{\mu}(\mathcal{C}_1)$ and $  \underset{n\rightarrow \infty}{\lim} n^{-\frac{d+2}{d+1}} \Gamma^n(\mathcal{C}_1) = \Gamma(\mathcal{C}_1)$. 
\end{prop}

In particular, this central limit theorem gives the generalization of the asymptotic rescaled mean and variance of $\XX^n_{\uu}$, as the definition (\ref{eq:limite_moyenne_covariance}) is almost identical to the rescaled asymptotic mean and variance of $\XX(\omega)$: we only restrict the integral on $\mathcal{C}_1$.

\begin{cor}
    Let $Z \in \mathcal{Z}(\mathcal{C}, n\kk)$ be a random zonotope drawn under $\Pn$, and $\XX^n_{\uu}$ the tangent point of the boundary of $Z$ to $H_{\uu}$. Then $\XX^n_{\uu}$ satisfies a central limit theorem in the sense that
$$ (\Gamma_{\uu}^n) ^{-1/2}\left(\XX^n_{\uu} -\boldsymbol{\mu}_{\uu}^n \right) \overset{(d)}{\underset{n\rightarrow \infty}{\longrightarrow}} \mathcal{N}(\boldsymbol{0}, I_d),$$
with $ \underset{n\rightarrow \infty}{\lim} \frac{1}{n} \boldsymbol{\mu}^n_{\uu} = \boldsymbol{\mu}_\uu$ and $  \underset{n\rightarrow \infty}{\lim} n^{-\frac{d+2}{d+1}} \Gamma_{\uu}^n = \Gamma_\uu$. 
\end{cor}

\begin{proof}
Proposition \ref{thm:CLT} is totally analogous to Proposition 4.1 in \cite{Barany-Bureaux:zonotopes} with the cone $\mathcal{C}_1$, except that $\kk$, which is the expected endpoint of a re-scaled zonotope drawn under $\Pn$ and which gives the vector $\aa$ in $\Pn$, does not belong to the cone $\mathcal{C}_1$. The proof relies on the Lyapunov condition that is a direct consequence of Lemma \ref{lem_Lyapu}. \medskip 

 We start with writing the expectation. After the expansion of the quotient as a series, the Fubini-Tonelli theorem yields:
$$ \boldsymbol{\mu}^n(\mathcal{C}_1) = \sum_{\xx \in \P_d \cap \mathcal{C}_1 } \frac{\xx e^{- \beta_n \aa \cdot \xx}}{1 - e^{- \beta_n \aa \cdot \xx}} = \sum_{i \geq 1} \sum_{\xx \in \P_d \cap \mathcal{C}_1 } \xx e^{- i \beta_n \aa \cdot \xx}.$$

For $i\leq 1/ \beta_n$, we use Proposition \ref{prop:equivalent_int_sum} to approximate the $i$-th summation over $\P_d \cap \mathcal{C}_1$ into a $d$-dimensional integral. Denoting $i_0$ the first integer such that $i_0 > 1/ \beta_n$,  there exists $a >0$ independent of $\xx$ such that 
\begin{align}\label{equ:summation_negligeable_CLT}
     \sum_{i > 1/\beta_n } \xx e^{- i \beta_n \aa \cdot \xx}  \leq  \xx e^{- i_0 \beta_n \aa \cdot \xx} \frac{a}{\beta_n}
\end{align}
$$ $$

Therefore, using again Proposition \ref{prop:equivalent_int_sum} on the summation over $\P_d \cap \mathcal{C}_1$ of the right-end term of \ref{equ:summation_negligeable_CLT}, we obtain that the terms with $i > 1/ \beta_n$ only contribute for $O(1/\beta_n)$. We obtain, as $n$ goes to $+ \infty$,

$$ \boldsymbol{\mu}^n(\mathcal{C}_1) = \frac{\zeta(d+1)}{\beta_n ^{d+1} \zeta(d)} \int_{\mathcal{C}_1} \xx e^{\aa \cdot \xx} d\xx + O\left( \frac{1}{\beta_n^{d+1}} \right).$$

After simplification of the prefactor of the integral, the change of variable $\xx = t \xx'$ with $t = \aa \cdot \xx$ gives

$$ \frac{1}{n} \boldsymbol{\mu}^n(\mathcal{C}_1)  = (d+1)! \int_{\mathcal{C}( \aa \leq 1) \cap \mathcal{C}_1} \xx d\xx +  O\left( n ^{\frac{-1}{d+1}} \right) .$$

The details of  the calculation of the asymptotic behavior of the variance are exactly the same, namely\\
\begin{align}\label{expres:variance_u}
    \Gamma^n(\mathcal{C}_1) = \sum_{\substack{\xx \in \P_d \cap \mathcal{C}_1}} \xx \xx^\intercal \frac{ e^{- \beta_n \aa \cdot \xx}}{\left(1 - e^{- \beta_n \aa \cdot \xx}\right)^2} = n^{\frac{d+2}{d+1}} \left(\frac{\zeta(d)}{\zeta(d+1)}\right)^{\frac{1}{d+1}} (d+1)! \int_{\mathcal{C}( \aa \leq 1) \cap \mathcal{C}_1} \xx \xx^\intercal d\xx + O\left(n \right).
\end{align}

The central limit theorem is obtained by ensuring that the Lyapunov ratio $L_{n,\mathcal{C}_1}$, defined just below, tends to 0 as $n$ grows large: 
$$ L_{n,\mathcal{C}_1} = \sum_{\substack{\xx \in \P_d \cap \mathcal{C}_1}} \mathbb{E} \left[ \left\lVert \Gamma^n(\mathcal{C}_1) ^{-1/2} \left(\omega(\xx) - \mathbb{E}[\omega(\xx)]\right) \xx  \right\rVert^3\right].$$

The Lyapunov ratio tends to 0 as it is bounded from above by the Lyapunov ratio of the endpoint of the zonotope $Z$ in the proof of \cite[Proposition 4.1]{Barany-Bureaux:zonotopes}. It is also a direct consequence of the following lemma. \medskip

\end{proof}

In order to state a local limit theorem for $\XX^n(\mathcal{C}_1)$ (Subsection \ref{subsection:LLT}), we compute the order of approximation of the Lyapunov ratio:

\begin{lem}\label{lem_Lyapu}
$$ L_{n,\mathcal{C}_1} = O\left(n ^{-\frac{d}{2(d+1)}}\right)$$
\end{lem}

\begin{proof}
We start by bounding the Lyapunov ratio with the operating norm of $\Gamma^n(\mathcal{C}_1)^{-1/2}$, denoted $||\Gamma^n(\mathcal{C}_1)^{-1/2} ||$:

$$\sum_{\substack{\xx \in \P_d \cap \mathcal{C}_1  }} \mathbb{E} \left[ \left\lVert \Gamma^n(\mathcal{C}_1) ^{-1/2} \left(\omega(\xx) - \mathbb{E}[\omega(\xx)]\right) \xx  \right\rVert^3\right] \leq ||\Gamma^n(\mathcal{C}_1) ^{-1/2} ||^3 \sum_{\substack{\xx \in \P_d  \cap \mathcal{C}_1 }} ||\xx||^3 \mathbb{E} \left[ \left| \left(\omega(\xx) - \mathbb{E}[\omega(\xx)]\right)   \right|^3\right]    $$

We respectively compute the second and fourth moments of $\overline{\omega}(\xx) = \omega(\xx) - \mathbb{E}[\omega(\xx)]$:

\begin{align*}
    \mathbb{E} \left[ \left|\overline{\omega}(\xx)   \right|^2\right] &= \frac{e^{-\beta_n \aa \cdot \xx}}{\left( 1 - e^{-\beta_n \aa \cdot \xx} \right)^2} \\
    \mathbb{E} \left[ \left|\overline{\omega}(\xx)   \right|^4 \right] &= \frac{ e^{-\beta_n \aa \cdot \xx} \left( 1 + 7 e^{-\beta_n \aa \cdot \xx} + e^{- 2\beta_n \aa \cdot \xx} \right)}{\left( 1 - e^{-\beta_n \aa \cdot \xx} \right)^4} \leq \frac{9 e^{-\beta_n \aa \cdot \xx}}{\left( 1 - e^{-\beta_n \aa \cdot \xx} \right)^4 }
\end{align*}

Therefore by applying the Cauchy-Schwarz inequality, we obtain

$$\mathbb{E} \left[ \left| \left(\omega(\xx) - \mathbb{E}[\omega(\xx)]\right)   \right|^3\right]  \leq \frac{3 e^{-\beta_n \aa \cdot \xx}}{\left( 1 - e^{-\beta_n \aa \cdot \xx} \right)^3 } $$

It is already known that $||\Gamma^n(\mathcal{C}_1)^{-1/2}||$ is of order $n^{-\frac{d+2}{2(d+1)}}$, therefore the end of the proof is obtained using the same arguments as in the computation of the asymptotic mean in the proof of Proposition \ref{thm:CLT}, namely

$$ \sum_{\substack{\xx \in \P_d \cap \mathcal{C}_1  }} ||\xx||^3 \frac{3 e^{-\beta_n \aa \cdot \xx}}{\left( 1 - e^{-\beta_n \aa \cdot \xx} \right)^3 }  = O\left( n^{\frac{d+3}{d+1}}\right).  $$
\end{proof}


We  extend the central limit theorem to the weak convergence of a $k$-tuple of rescaled tangent points to a Gaussian $k$-tuple. Let us introduce the rescaled point  $\boldsymbol{\Tilde{\XX}}^n_{\uu}$ of the boundary of the zonotope $Z$ tangent to the hyperplane $H_{\uu}$: 

\begin{align}\label{def:XX_tilde}
    \boldsymbol{\Tilde{\XX}}^n_{\uu} = n^{- \frac{d+2}{2(d+1)}} \left(  \XX^n_{\uu} - \mathbb{E}[ \XX^n_{\uu} ] \right)
\end{align}

\begin{cor}[Limit of a $k$-tuple of tangent points.]\label{prop:k-tuple-conv}
Let $Z$ be a zonotope drawn under $P_n$, and let $(\uu_1, ..., \uu_m)$, respectively $\left(\boldsymbol{\Tilde{\XX}}^n_{\uu_1}, ..., \boldsymbol{\Tilde{\XX}}^n_{\uu_m}\right)$, be an $m$-tuple of vectors of $\R^d\setminus\{\boldsymbol{0}\}$, respectively the rescaled $m$-tuple of the points of contact between $Z$ and the hyperplanes $H_{\uu_i}$, with $1 \leq i \leq m$. \medskip

Then 
$$\left(\boldsymbol{\Tilde{\XX}}^n_{\uu_1}, ..., \boldsymbol{\Tilde{\XX}}^n_{\uu_m}\right) \overset{\Pn}{\underset{n  \rightarrow + \infty}{\longrightarrow} } \left(\boldsymbol{\Tilde{\XX}}_{\uu_1}, ..., \boldsymbol{\Tilde{\XX}}_{\uu_m}\right)$$

where $\left(\boldsymbol{\Tilde{\XX}}_{\uu_1}, ..., \boldsymbol{\Tilde{\XX}}_{\uu_m}\right)$ is a centered Gaussian vector with covariance structure given by  
$$\text{Cov}(\boldsymbol{\Tilde{\XX}}_{\uu_i}, \boldsymbol{\Tilde{\XX}}_{\uu_j}) = \Gamma(\mathcal{C}(\uu_i \geq 0, \uu_j \geq 0)), \text{ for $1 \leq i,j \leq m$.} $$
\end{cor}

\begin{proof}
For $1 \leq i \leq k$, we write $\boldsymbol{\Tilde{\XX}}^n_{\uu_i}$ as the sum of all the generators contributing, that is:

$$\boldsymbol{\Tilde{\XX}}^n_{\uu_i} = n^{- \frac{d+2}{2(d+1)}} \sum_{\vv \in \P_d \cap \mathcal{C}(\uu_i \geq 0)}\left( \omega(\vv) - \mathbb{E}[\omega(\vv)] \right) \vv. $$

We introduce the sets $\left( \mathcal{C}_I \right)_{I \subset \llbracket 1, m \rrbracket}$ and $\left( \AA_I^n \right)_{I \subset \llbracket 1, m \rrbracket}$, respectively defined by, for $I \subset \llbracket 1, m \rrbracket$:

\begin{align}\label{equ:def_AA_I}
    \mathcal{C}_{I} = \underset{i\in I}{\bigcap} \mathcal{C}(\uu_i  \geq 0) \underset{j \notin I}{\bigcap} \mathcal{C}(\uu_j  < 0),\:\:\:\:\: \text{ and }\AA_{I}^n = n^{- \frac{d+2}{2(d+1)}} \sum_{\vv \in \P_d \cap \mathcal{C}_I}\left( \omega(\vv) - \mathbb{E}[\omega(\vv)] \right) \vv.
\end{align}

The cones $(\mathcal{C}_I)_{I \subset \llbracket1, m\rrbracket}$ are partitioning the cone $\mathcal{C}$, in particular they are mutually disjoint. As a consequence, the product form of $\Pn$ given in (\ref{eq:P_n_definition}) implies that the random variables $(\AA_I^n)$ are independent. Therefore, using Proposition \ref{thm:CLT},

$$ (\AA_I^n)_{I \subset \llbracket 1, m \rrbracket} \underset{n \rightarrow + \infty}{\overset{\Pn}{\longrightarrow}} \left( \AA_I \right)_{I \subset \llbracket 1, m \rrbracket}, $$

where $(\AA_I)_I$ are independent Gaussian variables of covariance $\boldsymbol{\Gamma} (\mathcal{C}_I)$. For $1 \leq i \leq m$, notice that we can write $\boldsymbol{\Tilde{\XX}}^n_{\uu_i}$ as a sum of $\AA^n_I$:

$$ \boldsymbol{\Tilde{\XX}}^n_{\uu_i} = \sum_{\substack{I \subset \llbracket 1, m \rrbracket \\ i\in I}} \AA_I^n.$$

Hence, for any $i, j$ between $1$ and $m$, the covariance of $\boldsymbol{\Tilde{\XX}}^n_{\uu_i}$ and $\boldsymbol{\Tilde{\XX}}^n_{\uu_j}$ is the variance of $\sum_{\substack{I \subset \llbracket 1, m \rrbracket \\ \{i, j \} \subset I}} \AA_I^n$. The weak convergence of the $\left(\boldsymbol{\Tilde{\XX}}^n_{\uu_i} \right)_i$ follows.\medskip

\end{proof}

\subsection{Local limit theorem}\label{subsection:LLT}
The central limit theorem is enough to obtain the limit shape of uniformly distributed zonotopes (\cite{Barany-Bureaux:zonotopes}), but not to ensure a Donsker theorem for the integral zonotopes. The local limit theorem below refines the approximation of the asymptotic behavior of the endpoint of generators in a $d$-dimensional subcone under $\Pn$. 

\begin{thm}[Local Limit Theorem]\label{thm:LLT}
Let $Z$ be a random integral zonotope drawn under the law $\Pn$. Let $\mathcal{C}_1 \subset \mathcal{C}$  be a $d$-dimensional cone with $0 \in \mathcal{C}_1$, and $\XX^n(\mathcal{C}_1)$ be the endpoint of the generators of $Z$ in $\mathcal{C}_1$. Then the random variable $\XX^n(\mathcal{C}_1)$ satisfies a local limit theorem of rate $n^{- \frac{d}{2(d+1)}}$. Formally:

$$\underset{n \rightarrow + \infty}{\lim\sup}\underset{\xx \in \Z_+^d}{\sup} n^{\frac{d}{2(d+1)}} \left|\Pn (\XX^n(\mathcal{C}_1) = \xx) - \frac{g_d \left( (\xx - \boldsymbol{\mu}^n(\mathcal{C}_1))^\top \Gamma^{n}(\mathcal{C}_1)^{-1} (\xx - \boldsymbol{\mu}^n(\mathcal{C}_1))\right)}{ \sqrt{ \det \Gamma^n(\mathcal{C}_1)}}  \right| < + \infty,   $$

where $g_d$ is the density of a standard normal $d$-dimensional variable. 
\end{thm}

This theorem is proved using the framework developed by J. Bureaux in \cite{Bureaux:LLT}. The idea of this framework is to use the inversion formula of the characteristic function on the probability $\Pn (\XX^n_{\uu} = \xx)$, and decompose the difference onto 3 different domains, which involves satisfying 3 different conditions (among which is Lemma \ref{lem_Lyapu}). Additionally to the previous notation, we denote $\sigma_{n }(\mathcal{C}_1)^2 $ the smallest eigenvalue of $\Gamma^n(\mathcal{C}_1)$.\\

\begin{lem}\label{lem_cov}
With the notation above, the inverse of the minimal eigenvalue of the covariance matrix satisfies $$\frac{1}{\sigma_{n }(\mathcal{C}_1) \sqrt{\det\Gamma^n(\mathcal{C}_1)}} = O\left(n^{- \frac{3(d+2)}{2(d+1)}}\right) $$ 
\end{lem}

\begin{proof}
This lemma directly comes from the asymptotic estimate of $\Gamma^n(\mathcal{C}_1)$. As seen in the central limit theorem, the covariance matrix estimate is 
$$ \Gamma^n(\mathcal{C}_1) = n^{\frac{d+2}{d+1}} \left(\frac{\zeta(d)}{\zeta(d+1)}\right)^{\frac{1}{d+1}} (d+1)! \int_{\mathcal{C}( \aa \leq 1, u \geq 0)} \xx^\intercal \xx d\xx + O\left(n \right).$$

Hence, with a diagonalization of $n^{- \frac{d+1}{d+2}} \Gamma^n(\mathcal{C}_1)$, we obtain $\sigma_{n}(\mathcal{C}_1)^2 \asymp n^{\frac{d+2}{d+1}}$, it follows that:

$$ \sigma_{n }(\mathcal{C}_1) \sqrt{\det(\Gamma^n(\mathcal{C}_1))} = O \left( n^{ \frac{3(d+2)}{2(d+1)}} \right).$$ 
\end{proof}

The last condition of the local limit theorem consists in bounding the characteristic function out of an ellipsoid denoted $\epsilon _{n,\uu}$ and defined as:
$$\epsilon _{n,\mathcal{C}_1} = \left\{ \tt \in \mathbb{R}^d \: : \: || \Gamma^{n}(\mathcal{C}_1)^{1/2} \tt || \leq \frac{1}{4 L_{n,\mathcal{C}_1}} \right\}$$

\begin{lem}\label{lem_characfunc}If the cone $\mathcal{C}_1$ has dimension 2 or more,
$$ \underset{\tt \in [- \pi, \pi ]^2 \setminus \epsilon _n }{\sup} \left| \mathbb{E}\left[e^{i \tt \cdot \XX^n_{\uu} }\right]\right| =O(n^{-1}) $$
\end{lem}

\begin{proof}
The outline of the proof is quite standard and can be found in \cite{Bureaux:LLT}. For any complex number $z$, the following inequality holds:

$$ \left| \frac{1 - |z|}{1-z} \right| \leq \exp \left( \Re (z) - |z| \right)$$

We apply it to the characteristic function:

$$\left|  \mathbb{E} \left[e^{i \tt \cdot \XX^n(\mathcal{C}_1) }\right] \right| = \prod_{\xx \in \P_d \cap \mathcal{C}_1} \left| \frac{1 - e^{- \beta_n \aa \cdot \xx}}{1 - e^{-(\beta_n \aa -i \tt) \cdot \xx}} \right| \leq \exp\left( \Re \left( \sum_{\xx \in \P_d \cap \mathcal{C}_1} e^{- (\beta_n \aa - i \tt) \cdot \xx} \right) -  \sum_{\xx \in \P_d \cap \mathcal{C}_1} e^{- \beta_n \aa  \cdot \xx} \right).  $$

This can be rewritten using the cosine as

\begin{align}\label{cosinus_function_carac}
     \left|  \mathbb{E} \left[e^{i \tt \cdot \XX^n(\mathcal{C}_1) }\right] \right| \leq  \exp\left(  \sum_{\xx \in \P_d \cap \mathcal{C}_1} e^{- \beta_n \aa \cdot \xx} (\cos(\tt \cdot \xx)  - 1 )   \right).
\end{align}

To bound the summation of the right-hand side of (\ref{cosinus_function_carac}), we construct a sequence of $\xx$ such that $\xx \cdot \tt$ is small enough to be well approximated. Using the diagonalization of $\Gamma^{n}(\mathcal{C}_1)^{1/2} $, there exists a positive constant $c_1$ such that $|| \Gamma^{n}(\mathcal{C}_1)^{1/2} \tt || \leq c_1 n^{\frac{d+2}{2(d+1)}} ||\tt||$. Furthermore, since $L_{n, \mathcal{C}_1} = O\left(n^{-\frac{d}{2(d+2)}}\right)$, for every $\tt \notin \epsilon_n$, there exists a second positive constant $c_2$ such that $ ||\Gamma^{n}(\mathcal{C}_1)^{ 1/2} \tt || \geq c_2  n^{ \frac{d}{2(d+1)}} $. We deduce that there is a constant $A > 0$ such that

$$ \underset{1 \leq i \leq d}{\max}(|t_i|) \geq  A n^{- \frac{1}{d+1}} $$

In the sequel, we denote $\textbf{e}_i$ the canonical standard basis vector of the $i^{\text{th}}$ coordinate. Using the symmetry, we may assume that $|t_1| \geq a n^{- \frac{1}{d+1}} $, which means $t_1 \in [- \pi, -A n^{-\frac{1}{d+1}}] \cup [A n^{-\frac{1}{d+1} }, \pi ]$. Notice that such $t_1$ with the condition $\textbf{e}_1 \notin \mathcal{C}_1^{\perp}$ exists because the characteristic function would not depend on $t_1$. The rest of the proof consists in finding 2 arithmetic sequences of primitive vectors, whose common differences are far enough from each other to ensure the convergence of the scalar product of one sequence with $\tt$ towards a polynomial limit. \medskip

$\mathcal{C}_1$ is at least of dimension 2, so we state that $\textbf{e}_2 \notin \mathcal{C}_1^{\perp}$ without loss of generality. Therefore, there exists in the interior of the cone $\mathcal{C}_1$ a primitive vector $\xx_1$ such that:

$$ \left\{
    \begin{array}{lll}
        \xx_1 \cdot \ee_1 = p, \xx_1 \cdot \ee_2 = q, \: \text{ with } \: p \wedge q =1 \text{ and }(p+1)\wedge q=1\\
        \xx_2 = \xx_1 + \ee_1 \in \P^d \cap \mathcal{C}(\uu\geq 0)  \\
        \xx_1 + 2 \ee_1 \in \mathcal{C}(\uu \geq 0)
    \end{array}
\right.   $$

The arithmetic sequences $(\xx_{1,i})_{i \geq 1}$ and $(\xx_{2,i})_{i \geq 1}$ defined by $\xx_{\alpha, i} = i q \xx_{\epsilon} + \ee_1$ (for $\alpha \in \{1,2\}$), are both sequences of primitive vectors, due to the coprimality of $p$ and $q$ on one side, and $p+1$ and $q$ on the other. The term $\cos(\tt \cdot \left(i q \xx_{\epsilon} + \ee_1\right))$ is periodic with respect to $i$, and its period is $\frac{2 \pi }{|\tt \cdot (q\xx_{\alpha})|}$. We compute a lower bound for the difference between the two periods $t_1$ and $t_2$ of respectively $(\xx_{1,i})$ and $(\xx_{2,i})$. We have \medskip

$$\frac{2 \pi}{\tt \cdot (q\xx_{1})} - \frac{2 \pi}{\tt \cdot (q\xx_{2})} = \frac{2 \pi(q t_1) }{q^2(\tt \cdot \xx_1)(\tt\cdot \xx_2) \xx_1} \geq A' \frac{ n^{-\frac{1}{d+1}} }{q^2||\xx_2||\times ||\xx_1||},$$

with $A'\geq 0$. Therefore we have a constant $A_{\xx_1} \geq 0$, depending only on the choice of $\xx_1$, such that at least one of these sequences has a period that differs from $1$ by $\frac{A_{\xx_1}}{2} n^{-\frac{1}{d+1}}$ or more. Similarly, both periods cannot be greater than $\frac{4 \pi}{A}n^{\frac{1}{d+1}}$ at the same time. For if we suppose that $|\tt \cdot \xx_1| \leq  \frac{A}{2} n^{\frac{1}{d+1}}$, then 

$$ |\tt \cdot \xx_2| \geq q|t_1| -  |\tt \cdot \xx_1|  \geq  \frac{A}{2} n^{\frac{1}{d+1}}.$$

Suppose $\xx_1$ verifies these conditions.  We can finally compute an upper bound for the argument of the exponential in \ref{cosinus_function_carac}:

$$  \sum_{\xx \in \P^d \cap \mathcal{C}_1} e^{- \beta_n \aa \cdot \xx} (\cos(\tt \cdot \xx)  - 1 )   \leq \sum_{i \geq 1} e^{- \beta_n \aa \cdot \xx_{1,i}} (\cos( \xx_{1,i} \cdot \tt )  - 1 ). $$

Denote $A_{\max} = \max(\frac{2}{A_\xx}, \frac{4 \pi}{A})$. The inequality $\cos(i \tt \cdot q\xx_1 + t_1) \leq \frac{1}{2}$ stands in the window of length three-quarters of the length of the period. The condition on the upper bound and the condition on the difference to 1 implies that the $k^{\text{th}}$ term of $(\xx_{1,i})_{i\geq 1}$ that verifies $\cos(\tt \cdot \xx_{1,i} ) \leq \frac{1}{2}$ is in the first $2k$ terms, for $2k \geq A_{\max}n^{\frac{1}{d+1}}$. This leads to:

$$\sum_{i \geq 1} e^{- \beta_n \aa \cdot \xx_{1,i}} (\cos( \xx_{1,i} \cdot \tt )  - 1 ) \leq  \left(- \frac{1}{2} \right) \sum_{i \geq A_{\max}n^{\frac{1}{d+1}}} e^{- \beta_n \aa \cdot \xx_{1,2i}} .   $$

Ultimately, a manipulation over the indicial notation gives:

$$ - \frac{1}{2} \sum_{i \geq A_{\max}n^{\frac{1}{d+1}}} e^{- \beta_n \aa \cdot \xx_{1,2i}} \leq - \frac{1}{2} \exp \left(- \beta_n \aa \cdot (A_{\max}n^{\frac{1}{d+1}} q\xx_{1} + \ee_1 )\right)  \sum_{i \geq 0} e^{- \beta_n i q \aa \cdot \xx_1} .   $$

We recall that
 $\beta_n =  \left(\frac{\zeta(d+1)}{\zeta(d) n }\right)^{ \frac{1}{d+1}}$, 
hence the first exponential asymptotically converges to a constant. A quick asymptotic analysis of the sum gives 

$$ \sum_{i \geq 0} e^{- \beta_n i q \aa \cdot \xx_1} = \frac{1}{1 - e^{- \beta_n q \aa \cdot \xx_1}} \asymp n^{\frac{1}{d+1}}. $$

Thus the biggest $\xx_1$ over all possible combinations of coordinates gives a constant $\gamma>0$, depending only on $\mathcal{C}_1$, such that 

$$\left|  \mathbb{E} \left[e^{i \tt \cdot \XX^n(\mathcal{C}_1)}\right] \right| \leq  \exp\left(- \gamma n^{\frac{1}{d+1}} \right),$$

which concludes the proof.

\end{proof}

\begin{proof}[Proof of the Theorem] 
Let $(a_n)_{n \in \N}$ be the sequence given by $a_n = n^{-\frac{d}{2(d+1)}}$. Then, with Lemma \ref{lem_Lyapu}, \ref{lem_cov}, and \ref{lem_characfunc}, the assumptions for Proposition 7.1 from \cite{Bureaux:LLT} are satisfied and there is a local limit theorem of rate $(a_n)$ for the variable $\XX^n_{\uu}$ under $P_n$. \medskip

\end{proof}

The following proposition gives the weak convergence of finite-dimensional distribution of the process of tangent points under $\text{Q}_{n\kk}$, leading to Theorem \ref{thm:1_CLT_under_Qn} and the Donsker theorem (Theorem \ref{thm:Donsker}). Recall that $\text{Q}_{n\kk}$ is the uniform distribution over integral zonotopes ending at $n\kk \in  \Z_+^d$. The connection between $\text{Q}_n$ and $\Pn$, for any $A \subset \Omega$:
$$\text{Q}_n [A] =  \Pn \left[A \:|\: \XX(\omega) = n \kk, \:\:\omega \in A\right] = \frac{\Pn \left[A \cap \{\XX(\omega) = n \kk,\:\:\omega \in A\} \right]}{\Pn\left[\XX(\omega) = n \kk,\:\:\omega \in A\right]}.$$

\begin{prop}[Weak convergence of finite-dimensional marginals]\label{prop:finite-dim-distri}
For $\kk \in \text{int }\mathcal{C} \cap \Z^d$, let $(\uu_1, ..., \uu_m)$ an $m$-tuple of $\R^d$, let $Z$ be a random zonotope drawn under $\text{Q}_{n\kk}$, and $\left(\boldsymbol{\Tilde{\XX}}^n_{\uu_1},..., \boldsymbol{\Tilde{\XX}}^n_{\uu_m}\right)$ be the rescaled position of the tangent points of $Z$ and $H_{\uu_1}$, ..., $H_{\uu_m}$. Then there is an independent family of Gaussian centered random variables $\left(\GG_{I}\right)_{I \subset \llbracket 1, m\rrbracket}$ such that \\
\begin{align}
    \left(\boldsymbol{\Tilde{\XX}}^n_{\uu_1},..., \boldsymbol{\Tilde{\XX}}^n_{\uu_m}\right) \overset{(d)}{\underset{n \rightarrow + \infty}{ \longrightarrow}} (\NN_{\uu_1},..., \NN_{\uu_m}),
\end{align}

where $\NN_{\uu_i} = \sum_{\substack{I \subset\llbracket1, m \rrbracket \\i \in I}}\GG_{I}$ and 
$${\normalfont \text{Cov}}(\GG_{I}) = {\footnotesize \left( \Gamma(\underset{i\in I}{\cap} \mathcal{C}(\uu_i  \geq 0) \underset{j \notin I}{\cap} \mathcal{C}(\uu_j  < 0))^{-1} + \Gamma(\underset{i\in I}{\cap} \mathcal{C}(\uu_i  < 0) \underset{j \notin I}{\cap} \mathcal{C}(\uu_j  \geq 0))^{-1} \right)^{-1} }$$ When $\mathcal{C}_{u_i}$ is $d$-dimensional, $\NN_{\uu_i}$ is a centered Gaussian variable of variance $\left( \Gamma_{\uu_1}^{\:\:-1} + \Gamma_{-\uu_i} ^{\:\:-1} \right)^{-1}$.
\end{prop}  

\begin{proof} 
We still denote $\omega$ as the function of multiplicities of $Z$. Starting with $\uu_1, ...,\uu_m$ vectors of $\R^d$, we denote $\uu_{m+1}$ a vector in the dual cone of $\mathcal{C}$ (that is a vector $\vv$ such that $\vv \cdot \xx >0 $ for all $\xx \in \mathcal{C}$), hence $\boldsymbol{\Tilde{\XX}}^n_{\uu_{m+1}}$ is the rescaled endpoint of the zonotope $Z$. \medskip

Using the same argument as in Proposition \ref{prop:k-tuple-conv}, the probability measure $\Pn$ is constructed as a product of geometric distributions of the primitive vectors in $\mathcal{C}$, and the occurrences of the primitive vectors (that is the set $\{\omega(\xx)\}_{\xx}$) are mutually independent. Therefore we introduce again the family of variables $\left(\AA_{I}\right)_{I \subset \llbracket1 , m\rrbracket}$ defined in (\ref{equ:def_AA_I}).\medskip


$\left( \AA_I\right)_{I \subset \llbracket1, m\rrbracket }$ denotes the vertices of $Z$ at the end of a path starting at the origin and composed of all generators contributing to the elements of I, but not contributing to the others, after recentering and renormalizing  by $n^{-{(d+2)}/{(2d+2)}}$. The family of cones $\left(\mathcal{C}_I\right)_{I \subset \llbracket1, m\rrbracket}$ is mutually disjoint, hence the variables $\left(\AA_{I}\right)_{I \subset \llbracket1, m\rrbracket}$ are mutually independent under $\Pn$. \medskip

The family $\left(\boldsymbol{\Tilde{\XX}}^n_{\uu_i} \right)$ is generated by the family $\left( \AA_I \right) _{I \subset \llbracket 1, m \rrbracket, \: I \neq \emptyset}$, as $\boldsymbol{\Tilde{\XX}}^n_{\uu_i} = \sum_{I \subset \llbracket 1, m \rrbracket, i \in I} \AA_I$. Under $\Pn$, the probability for the $(2^{m}-1)$-tuple of variables $\AA_I$, with $I \neq \emptyset$ to be equal to $\left(\xx_{I}\right)_{I \subset \llbracket1, m\rrbracket, \: I \neq \emptyset}$ is, based on Bayes' theorem:

\begin{align}\label{equ:Baye_AA_I}
    \text{Q}_n\left[\left(\AA_I\right)_{I \neq \emptyset} = \left(\xx_I \right)_{I \neq \emptyset } \right] = \frac{\Pn \left[ \left(\AA_I\right)_I = \left(\xx_I \right)_I  \cap \boldsymbol{\Tilde{\XX}}^n_{\uu_{m+1}} = 0 \right]}{\Pn\left[\boldsymbol{\Tilde{\XX}}^n_{\uu_{m+1}} = 0\right]}.
\end{align}

Under the distribution $\text{Q}_n$, the condition that the endpoint of the zonotope $Z$ is at $n\kk$ can be considered as a condition on $\AA_{ \emptyset}$, that is $\AA_{\emptyset} = \xx_{\emptyset}$ with $\xx_{\emptyset} = -\sum_{I \subset \llbracket1, m\rrbracket, \: I \neq \emptyset } \xx_I$. Therefore, denoting $\mathcal{P}(E)$ the set of all subsets of $E$, it follows that

\begin{align}\label{equ:indep_AA_I}
    \Pn \left[ \left(\AA_I\right)_{I \subset \llbracket1, m\rrbracket, \: I \neq \emptyset  } = \left(\xx_I \right)_{I \subset \llbracket1, m\rrbracket, \: I \neq \emptyset  }  \cap \boldsymbol{\Tilde{\XX}}^n_{\uu_{m+1}} = 0 \right] &= \Pn \left[ \AA_I = \xx_I,\:\: I \in \mathcal{P}\left(\llbracket 1, m \rrbracket \right) \right] \nonumber\\
    &= \prod_{I\in \mathcal{P}\left(\llbracket 1, m \rrbracket \right)} \Pn (\AA_I = \xx_I).
\end{align}



All the variables $\AA_I$ satisfy Theorem \ref{thm:LLT} with mean $0$ and covariance $\Gamma({\mathcal{C}_I})$, and so does $\boldsymbol{\Tilde{\XX}}^n_{\uu_{m+1}}$ with covariance $\Gamma_{\mathcal{C}}$. Hence, (\ref{equ:Baye_AA_I}), (\ref{equ:indep_AA_I}), and the local limit theorem lead to

\begin{align}\label{eq:LLT_pont}
    \underset{\left(\xx_I \right)_{\underset{I \neq \emptyset}{I \subset \llbracket1, m\rrbracket }} \in \left(\R^d\right)^{2^m-1} }{\sup} \left|
    \text{Q}_n\left[\left(\AA_I\right)_{I\neq \emptyset} 
    = \left(\xx_I \right)_{I \neq \emptyset} \right] - \frac{\prod_{I \subset \llbracket 1, m \rrbracket} \frac{g_d \left( \xx_I^\top (\Gamma({\mathcal{C}_I}))^{-1} \xx_I \right)}{ \det\Gamma({\mathcal{C}_I})^{1/2}} }
    {\frac{1}{\left(2 \pi \right)^{d/2} \det\Gamma({\mathcal{C}})^{1/2}}}\right| \underset{n \rightarrow + \infty}{\longrightarrow} 0,
\end{align}

where $g_d$ is the density of the $d$-dimensional standard Gaussian variable. $\left(\AA_I\right)_{I \subset \llbracket 1,m \rrbracket, \: I \neq \emptyset}$ satisfies a local limit theorem of rate 1 to a Gaussian family of variables named $\left(\GG_{I} \right)_{I \subset \llbracket 1,m \rrbracket, \: I \neq \emptyset}$. The weak convergence follows. The covariance of the $(\GG_I)$ is given by inverting the matrix of the quadratic form in the exponential. \medskip

Since $\boldsymbol{\Tilde{\XX}}^n_{\uu_i} = \sum_{I \subset \llbracket 1, m \rrbracket, i \in I} \XX_I$, $\left(\boldsymbol{\Tilde{\XX}}^n_{\uu_i}\right)_{1 \leq i \leq m}$ weakly converges to $\left( \sum_{\substack{I_i \subset \llbracket 1, m \rrbracket \\ i \in I_i}}\GG_{I_i} \right)_{1 \leq i \leq m}$. When considering this limit with only one variable $\boldsymbol{\Tilde{\XX}}^n_{\uu_i}$, we deduce that it weakly converges to the Gaussian random variable $\NN_{\uu_i}$ with mean $0$ and covariance

\begin{align}
    \text{Cov}\left(\NN_{\uu_i}\right) =\left( \Gamma(\mathcal{C}(\uu_i \geq 0))^{-1} + \Gamma(\mathcal{C}\setminus \mathcal{C}(\uu_i \geq 0)) ^{-1} \right)^{-1} = \left( \Gamma_{\uu_i}^{-1} + \Gamma_{-\uu_i} ^{-1} \right)^{-1}.
\end{align}

\end{proof}

\begin{proof}[Proof of Theorem \ref{thm:1_CLT_under_Qn}]
    Theorem \ref{thm:1_CLT_under_Qn} is a direct consequence of Proposition \ref{prop:finite-dim-distri}, taking only one vector $\uu$. 
\end{proof}

%% file: Parties/Development/3donsker.tex
This part is dedicated to proving the weak convergence of the rescaled process of the fluctuations around the limit shape of a random zonogon to a Brownian bridge, in the same vein as Donsker's theorem. This will be used in the next section to prove Theorem \ref{thm:brownian_polygone} on random polygons.

Recall that the results of the previous section were true in any dimension. On the other hand, extending the functional results of this section to higher dimension would require more involved tightness estimates, such as those appearing in Theorem 4 of \cite{bickel:multidim_tightness}. We will not handle this question here. In any case, in dimension $\geq 3$, the only results we could obtain using our method would be valid for zonotopes and not for polytopes, whereas in dimension 2, the study of polygons follows readily from the case of zonogons.  \\

So in this section, $d =2$ and $\mathcal{C}= \R_+^2$. The boundary of a zonogon $Z$ of endpoint $n\kk$ is divided into two polygonal lines, identical up to central symmetry, denoted $Z^+$, resp. $Z^-$, for the upper polygonal line of $Z$, resp. for the lower line. The tangent point of $Z$ to the plane $H_{\uu}$ belongs to $Z^+$ if $\uu \in \R_{-} \times \R_{+}$, to $Z^{-}$ if $\uu \in \R_{+} \times \R_{-}$, and it is $n\kk$ if $\uu \in \R_{+}^2$ and $0$ if $\uu \in \R_{-}^2$. We restrict the study to $\Z^+$ and $\uu \in \R_{-} \times \R_{+}$, and we define the process $\left(\BB^n_{t}, 0 \leq t \leq 1 \right)$ of fluctuations of tangent points of $Z^+$ away from their mean position, where $\BB^n_{t} = \boldsymbol{\Tilde{\XX}}^n_{(t-1, t)}$. \medskip







\begin{figure}[H]
    \centering
    \begin{overpic}[width = 10cm]{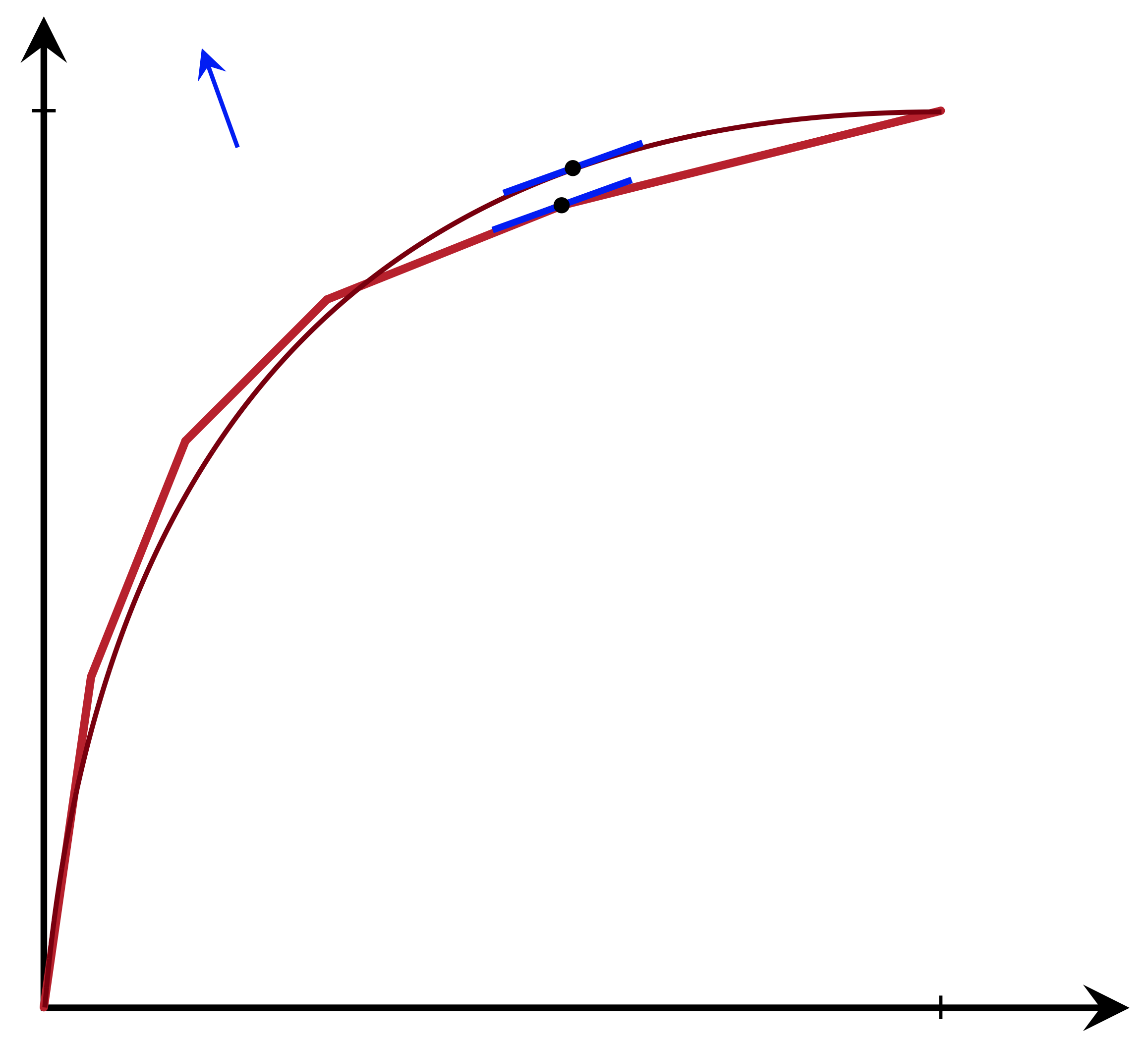}
        \put(80.5,-1.3){$n$}
        \put(-1,80){$n$}
        \put(12,74){\color{blue}$\uu = (t-1, t)$}
        \put(46,68){$\XX_{\uu}^n$}
        \put(45,79){$\mathbb{E}\left[ \XX^n_\uu \right]$}
    \end{overpic} %
    \caption{Upper polygonal line of a zonotope and its limit shape in the square $[0,n]^2$}%
    \label{fig:ecartalacourbe}%
\end{figure}
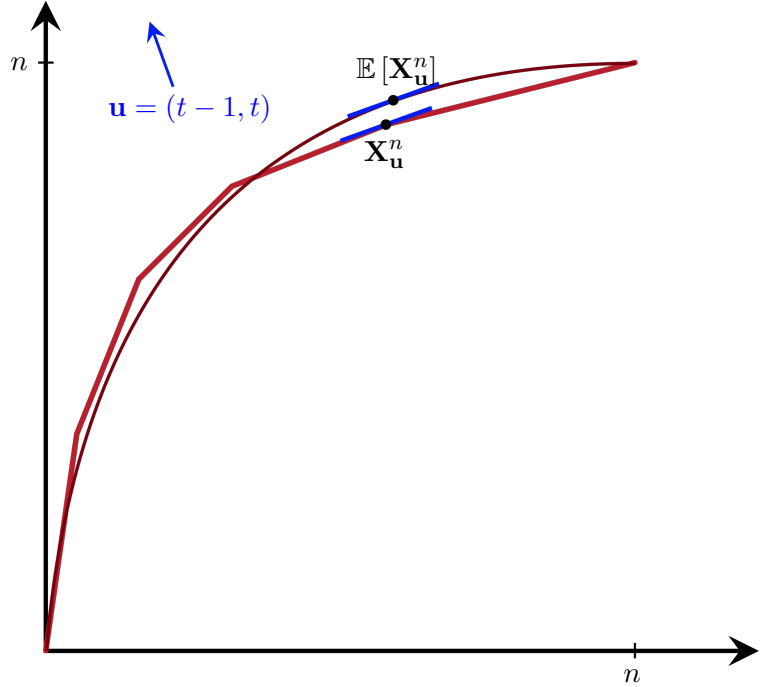

\begin{thm}\label{thm:Donsker}
Let $Z$ be a random, uniform zonogon in $\mathcal{Z}( \R_+^2, (n,n))$. For $t\in[0,1]$, put $\BB^n_{t} = \boldsymbol{\Tilde{\XX}}^n_{(t-1, t)}$
with $\boldsymbol{\Tilde{\XX}}^n_{(t-1, t)}$ as defined in \eqref{def:XX_tilde}. Then we have the convergence of processes 
$$\left(\BB^n_{ t}, 0 \leq t \leq 1 \right)\stackrel{(d)}{\to} (\PP_t)$$ on the space $D[0, 1]$  of c\`adl\`ag functions on $[0,1]$ equipped
with the Skorokhod topology, where  $(\PP_t)$ is a Gaussian processus with covariance matrix given by \eqref{covariance}. Alternatively,   $(\PP_t)$ can be written as $\PP_t=  P(t)D(t)\BB_t $, where $D(t)$ is a non-degenerate diagonal matrix, $P(t)$ is an orthogonal matrix and $\left( \BB_t\right) _{t\in [0,1]}$ is the 2-dimensional standard Brownian bridge.  
\end{thm}

\begin{proof}
    The weak convergence of the finite-dimensional distributions to those of $(\PP_t)$ is given by Proposition \ref{prop:finite-dim-distri}, for $d= 2$. For $0 < t < 1$, the covariance of $\PP_t$ is $(\Gamma_{(t-1, t)}^{-1} + \Gamma_{(1-t, -t)} ^{-1})^{-1}$, that is

\begin{align}\label{covariance}
    \text{Cov}(\PP_t) = \left(\frac{\zeta(2)}{\zeta(3)}\right)^{1/3}  \begin{pmatrix}
-2 \left((4t^2-2t+1\right)t^3(t-1) & \left( 8t^2 - 8t +3 \right) (t-1)^2 t^2 \\
\left( 8t^2 - 8t +3 \right) (t-1)^2 t^2 & -2\left(4t^2 - 6t +3 \right)(t-1)^3 t
\end{pmatrix}.
\end{align}

Using the spectral theorem, we compute the orthogonal matrix $P(t)$ such that $P(t)\text{Cov}(\PP_t) P(t)^{\intercal}$ is diagonal. We define the two polynomials $f(x) = - 8 x^2 - 8 x -1/2 $ and $g(x) = 64 x^4 + 128 x^3 +69 x^2 + 19/2 x + 1/16$, and we obtain

$$ P(t) \text{Cov}(\PP_t) P(t)^{\intercal} =  t(1-t)\begin{pmatrix}
Q_{-}(t)  & 0\\
0& Q_{+}(t)
\end{pmatrix}, $$

where $Q_{\pm}(t) = f\left( (t-\frac{1}{2})^2 \right) \pm \sqrt{g\left( (t-\frac{1}{2})^2 \right)}$. 

We only need to prove that $\left(\BB^n_{ t}, 0 \leq t \leq 1 \right) $ is tight, which is a consequence of the following proposition, using a result of Billingsley (Theorem 13.5 in \cite{billingsley:converg}).\medskip





\end{proof}

Remarkably, $Q_{-}(t) Q_{+}(t) = 3 t(1-t)$. These functions are both symmetrical to $t = 1/2$, and $Q_1$ cancels out at $0$ and $1$ while $Q_2$ cancels out at $1/2$. We don't have any interpretation yet for these terms.  We shall not give the explicit formula of the orthogonal matrix $P(t)$ here but we display the asymptotic behavior at $t=0$ and $t=1$:

\begin{align}\label{orthogonal}
    P(t) &=  \begin{pmatrix}
-1 + \frac{1}{8}t^2  & \frac{1}{2}t + \frac{1}{6}t^2\\
\frac{1}{2}t + \frac{1}{6}t^2 &  1 - \frac{1}{8}t^2
\end{pmatrix} + o(t^2), \\
P(t) &= 
 \begin{pmatrix}
 \frac{1}{2}(t-1) + \frac{1}{6}(t-1)^2 & -1 + \frac{1}{8}(t-1)^2 \\
1 - \frac{1}{8}(t-1)^2 &  \frac{1}{2}(t-1) + \frac{1}{6}(t-1)^2
\end{pmatrix} + o\left((t-1)^2\right).\nonumber
\end{align}\medskip

\begin{prop}[Tightness]\label{prop:tightness}
For $ 0 \leq r < s < t \leq 1$, and $\alpha >\frac{1}{2}$, and $\beta > 0$

$$\mathbb{E}\left[ \left|\left|\BB^n_s - \BB^n_r \right|\right|_1^2 \left|\left| \BB^n_t-\BB^n_s \right|\right|_1^2 \right] \leq \left| F(t) - F(r)\right|^{2 \alpha}, $$
with $F$ is a non-decreasing, continuous function on $[0,1]$
\end{prop}

\begin{proof}
Let $0 \leq r \leq s \leq t\leq 1$, and . It is sufficient to prove that there exists a constant $C >0$ such that

$$ \mathbb{E}\left[ \left|\left|\BB^n_s - \BB^n_r \right|\right|_1^2 \left|\left| \BB^n_t-\BB^n_s \right|\right|_1^2  \right] \leq C \left(t^3 - r^3\right)^2 .$$

We denote $\mathcal{E}(r,s,t)$ the left term of this inequality. Since $\{\BB^n_t\}$ is a variable on a random zonotope that end at $(n,n)$, the first step is to broaden to random zonotopes drawn under $\Pn$. As we are in 2 dimensions, we will write in the following $\XX^n_{t}$ instead of $\XX^n_{(1-t, t)}$, referring to the point of the boundary of $Z$ that is tangent to the hyperplane $H_{(1-t, t)}$. We have

$$\mathcal{E}(r,s,t) = \frac{1}{n^{8/3}} \mathbb{E}\left[ \left|\left|\XX^n_s - \XX^n_r - \mathbb{E}\left[\XX^n_s - \XX^n_r \right]  \right|\right|_1^2 \left|\left| \XX^n_t-\XX^n_s - \mathbb{E}\left[\XX^n_t - \XX^n_s \right] \right|\right|_1^2  \:\: \left| Z \in \mathcal{Z} \left(\mathbb{N}^2, (n,n) \right) \right.  \right]. $$

Denoting $\mathcal{C}_{r,s}$ (resp. $\mathcal{C}_{s,t}$)  the cone of vectors contributing to $\XX^n_s - \XX^n_r$ (resp. $\XX^n_t - \XX^n_s$). The first cone formally is $\left\{ \xx \in \R, \xx \cdot (r-1, r) \leq 0 \leq \xx \cdot (s-1, s) \right\}$. We can write the right term of the equation above as:

$$\mathcal{E}(r,s,t) = \mathbb{E}\left[ \left|\left| \boldsymbol{\Tilde{\XX}}^n(\mathcal{C}_{r,s})  \right|\right|_1^2 \left|\left| \boldsymbol{\Tilde{\XX}}^n(\mathcal{C}_{s,t}) \right|\right|_1^2  \:\: \left| Z \in \mathcal{Z} \left(\mathbb{N}^2, (n,n) \right) \right.  \right]. $$

After writing each variable as a sum over the generators and expanding the product of the norms, one can notice that this expectancy will result in a sum over quadruplets primitive generators $\{ \xx_1, \xx_2 \xx_3, \xx_4 \}$. Given such a quadruplet $\{ \xx_1, \xx_2, \xx_3, \xx_4 \}$, the probability for the generators to respectively occur $k_1, k_2, k_3$, and $k_4$ is

$$ \text{Q}_n(k_1, k_2, k_3, k_4) =\prod_{i=1}^4 \text{P}_n \left( \omega(\xx_i) = k_i \right) \frac{\text{P}_n \left(  \XX_1^n \setminus \{\xx_1,\xx_2,\xx_3,\xx_4\} = (n,n) - \sum_{i=1}^4 k_i \xx_i \right)}{\text{P}_n \left( \XX_1^n = (n,n)\right)} $$

Yet the local limit theorem ensures that there exist $n_1$ large enough such that for $n> n_1$, for any $\{ \xx_1, \xx_2, \xx_3, \xx_4 \}$, the ratio is bounded : 

$$ \frac{\text{P}_n \left(  \XX_1^n \setminus \{\xx_1,\xx_2,\xx_3,\xx_4\} = (n,n) - \sum_{i=1}^4 k_i \xx_i \right)}{\text{P}_n \left( \XX_1^n = (n,n)\right)} \leq 2 $$

Therefore, it remains

$$  \mathcal{E}(r,s,t) \leq 2\: \mathbb{E}\left[ \left|\left|\boldsymbol{\Tilde{\XX}}^n(\mathcal{C}_{r,s}) \right|\right|_1^2 \left|\left| \boldsymbol{\Tilde{\XX}}^n(\mathcal{C}_{s,t}) \right|\right|_1^2 \right] = 2\: \mathbb{E}\left[ \left|\left|\boldsymbol{\Tilde{\XX}}^n(\mathcal{C}_{r,s}) \right|\right|_1^2 \right] \mathbb{E}\left[\left|\left| \boldsymbol{\Tilde{\XX}}^n(\mathcal{C}_{s,t}) \right|\right|_1^2 \right]$$

 This expectancy can be split into two parts due to independence. In order to expand the 1-norm, let $a$ and $b$ be real numbers such that $\boldsymbol{\Tilde{\XX}}^n(\mathcal{C}_{r,s}) = (a,b)$. Therefore, $||\boldsymbol{\Tilde{\XX}}^n(\mathcal{C}_{r,s})||_1^2 = a^2 + b^2 + 2 |ab| \leq 2 a^2 + 2 b^2$. Hence, we have, writing $\xx = (x_1, x_2)$: 

 $$  \mathbb{E}\left[ \left|\left|\boldsymbol{\Tilde{\XX}}^n(\mathcal{C}_{r,s}) \right|\right|_1^2 \right]\leq \frac{2}{n^{4/3}} \:\mathbb{E}\left[ \left(\sum_{\xx \in \mathbb{P}^2 \cap \mathcal{C}_{r,s} } x_1(\omega(\xx) - \mathbb{E}[\omega(\xx)]) \right)^2 \right]
+\frac{2}{n^{4/3}}\:\mathbb{E}\left[ \left(\sum_{\xx \in \mathbb{P}^2 \cap \mathcal{C}_{r,s} } x_2(\omega(\xx) - \mathbb{E}[\omega(\xx)]) \right)^2 
\right] $$
 
 The 2 terms on the right-hand side are the diagonal terms of $\Gamma^n(\mathcal{C}_{r,s})$. Using the bounding given by Corollary \ref{cor:homogeneous_bound_on_cone} in the same way we previously explicitly calculated in the proof in Proposition \ref{thm:CLT}, we have

$$\mathbb{E}\left[ \left(\sum_{\xx \in \mathbb{P}^2 \cap \mathcal{C}_{r,s} } x_1(\omega(\xx) - \mathbb{E}[\omega(\xx)]) \right)^2 \right] \leq 6 n^{4/3} \left(\frac{\zeta(2)}{\zeta(3)}\right)^{\frac{1}{3}} \int_{\substack{\P_2 \cap \mathcal{C}_{r,s} \\ \xx \cdot (1,1) \leq 1}} x_1^2 d\xx + \epsilon(n)n, $$

with $\epsilon(n)$ tending to $0$ as $n$ grows to $+\infty$, and $\epsilon(n)$ is independent of $\mathcal{C}_{r,s}$ Therefore for $n_2$ large enough, and $n\geq n_2$, there exists a constant $A >0$ such that \medskip
 
$$ \mathbb{E}\left[ \left|\left|\boldsymbol{\Tilde{\XX}}^n(\mathcal{C}_{r,s}) \right|\right|_1^2 \right] \leq  A (s^3 - r^3) $$
 
 Finally, the inequality $(s^3 - r^3)(t^3 - s^3) \leq (t^3 - r^3)^2$ concludes the proof.

\end{proof}

In the proof of the convergence of the variations of a polygon uniformly drawn in a square, we need to extend Theorem \ref{thm:Donsker} to zonogons ending at $(n + r^{2/3}, n+ sn^{2/3})$. This is handled by the following proposition. 
\begin{prop}
Let $r, s\in\R$ and let $Z$ be a random, uniform zonogon in $\mathcal{Z}\left( \R_+^2, \left(\lfloor n + rn^{2/3} \rfloor\right.\right. \left.\left., \lfloor n + s n^{2/3} \rfloor \right)\right)$. Let $\BB_t^n$ be, as in Theorem 3, 
the rescaled position of the point of $Z$ tangent to the hyperplane normal to $(t-1, t)$. Then we have the convergence of processes
$$\left(\BB^n_{ t}, 0 \leq t \leq 1 \right)\stackrel{(d)}{\to} \PP_t + \boldsymbol{\nu}_{r, s}(t)$$
in the space $D[0,1]$ of c\`adl\`ag functions on $[0,1]$ equipped
with the Skorokhod topology,  where $(\PP_t)$ is the same process as in Theorem \ref{thm:Donsker} and $\boldsymbol{\nu}_{r, s}(t)$ is a drift term given by \eqref{nu}.\\
Denoting $a = r-s$ and $b=r -2s$, $\left(\boldsymbol{\nu}_{r, s}(t) \right)_{0\leq t \leq 1}$ is a parametric cubic curve starting at $(0,0)$, ending at $(r,s)$ and satisfying the equation:
 \begin{align}\label{eq:equation_nu}
      4a(X + Y)^3 = 28b^3 X +7 b^2(X + Y)^2 + 54ab X(X+Y) + 27a^2 X^2 
 \end{align}
\end{prop}

\begin{figure}[H]
        \centering
        \subfloat[\centering $r=0.1$, $s=1$  ] {{ \includegraphics[width=5cm]{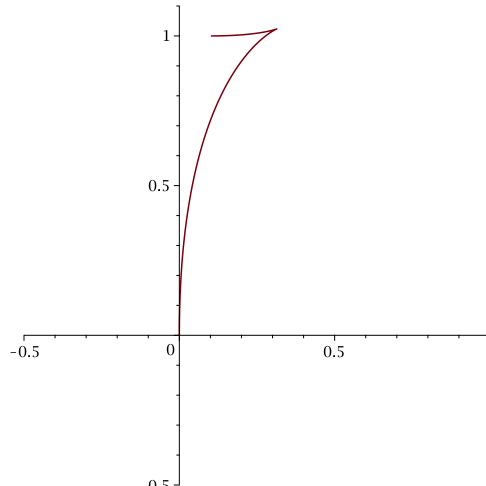} }}%
        \qquad
        \subfloat[\centering $r=1$, $s=1$ ] {{ \includegraphics[width=5cm]{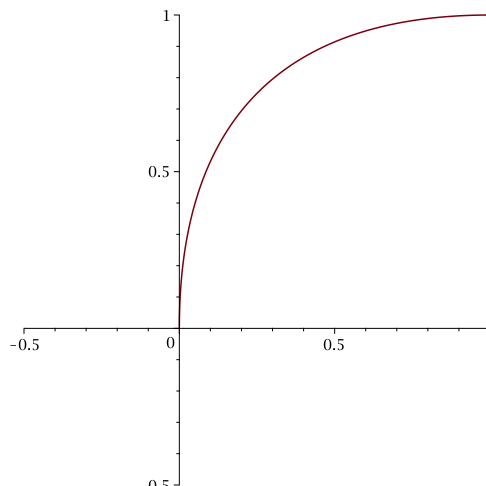} }}%
        \qquad
        \subfloat[\centering $r=1$, $s=0.5$  ] {{ \includegraphics[width=5cm]{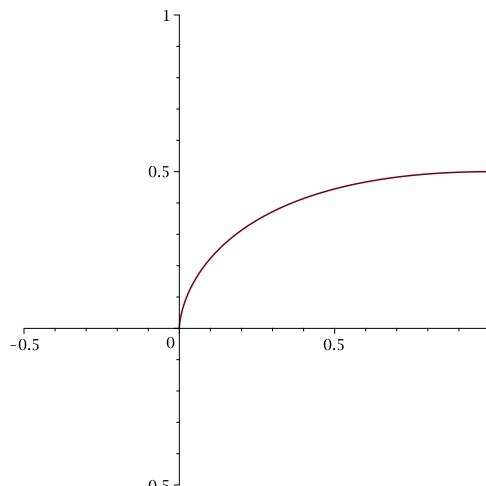} }}%
        \qquad
        \subfloat[\centering $r=1$, $s=0.2$  ] {{ \includegraphics[width=5cm]{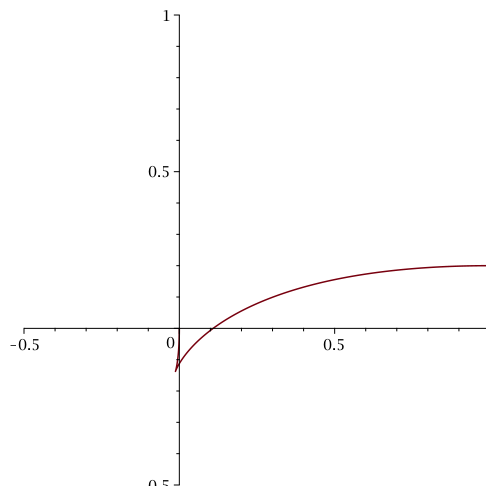} }}%
        \qquad
        \subfloat[\centering $r=1$, $s=-0.4$ ] {{ \includegraphics[width=5cm]{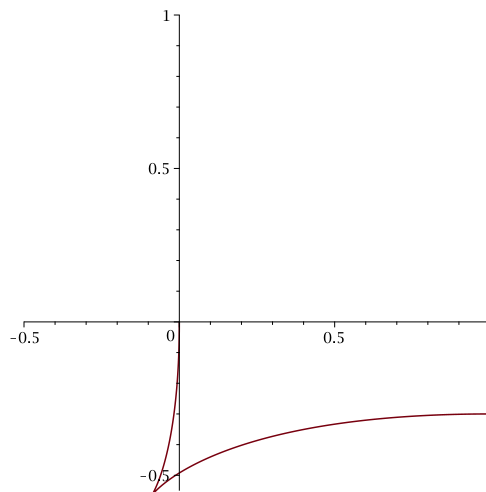} }}%
        \caption{The curve $\boldsymbol{\nu}_{r, s}(t)$ for $0 \leq t \leq 1$, for different values of $(r,s)$.}%
        \label{fig:nu_rs}%
    \end{figure}

\begin{proof}
     Both finite-dimensional distribution and tightness proofs are totally analogous to the proofs of Proposition \ref{prop:finite-dim-distri} and \ref{prop:tightness}. Given $\Gamma_{(t-1, t)}$ and $\Gamma_{-(t-1, t)}$ the covariances of $\boldsymbol{\Tilde{\XX}}_{(t-1, t)}$ and $\boldsymbol{\Tilde{\XX}}_{-(t-1, t)}$, the local limit theorem gives the following drift $\boldsymbol{\nu}_{r,s}(t)$:

    \begin{align}\label{nu}
         \boldsymbol{\nu}_{r, s}(t) &= \left(\Gamma_{(t-1,t)}^{\:-1} + \Gamma_{-(t-1, t)}^{\: -1} \right)^{-1} \Gamma_{-(t-1, t)}^{\:-1}  \begin{pmatrix} r \\ s \end{pmatrix} \nonumber\\
         &= \left(\begin{matrix}
                t^2(2t - 1) & -2 t^2(t-1)\\
                -2t(t-1)^2  & t(2t^2 -5 t + 4)
            \end{matrix}\right)
            \begin{pmatrix} r \\ s \end{pmatrix} 
    \end{align}

    All that remains is to analyze to find the equation of the curve $\boldsymbol{\nu}_{r, s}$. Notice that the coordinate of $\boldsymbol{\nu}_{r, s}$, $X(\boldsymbol{\nu}_{r, s})$ and $Y(\boldsymbol{\nu}_{r, s})$ are polynomials of $t$ of degree 3, but their sum is a polynomial of degree 2, that is $X(\boldsymbol{\nu}_{r, s})(t) +Y(\boldsymbol{\nu}_{r, s})(t) = 3at^2 -2 bt$ denoting $a = r-s$ and $b = r-2s$. \medskip

    In order to get an equation of the curve, one can resolve the polynomial of degree 2 and inject the solution in $X(\boldsymbol{\nu}_{r, s})(t)$. After simplification, we obtain the equation (\ref{eq:equation_nu}) 

    $$ 4a(X + Y)^3 = 28b^3 X +7 b^2(X + Y)^2 + 54ab X(X+Y) + 27a^2 X^2$$
    which is verified by $(X(\boldsymbol{\nu}_{r, s}),Y(\boldsymbol{\nu}_{r, s})$. When $r = s$, the curve is the parabola that is the limit shape of a uniform integral zonogon ending at $(r,s)$; otherwise, $\boldsymbol{\nu}_{r, s}(t)$ is a cubic curve. There are 2 different shapes, depending on $\frac{r}{s}$ belonging to $[1/2, 2]$ or not. \medskip The derivatives of $X(\boldsymbol{\nu}_{r, s})$ and $Y(\boldsymbol{\nu}_{r, s})$ with respect to $t$ are:
    $$ X(\boldsymbol{\nu}_{r, s})' (t) = 6 at^2 -2 b t, \:\:\: \text{ and }  Y(\boldsymbol{\nu}_{r, s})' (t) =- 6 at^2 + (6a + 2b)t -2b   $$
    
    If $\frac{r}{s} \in [1/2, 2]$, there is a cusp of multiplicity $2$ at $t_0= \frac{b}{3a} = \frac{r-2s}{3(r-s)}$. See Figure \ref{fig:nu_rs} for some examples of $(r,s)$.
\end{proof}

%% file: Parties/Development/5polygons.tex
The aim of this section is to prove Theorem 2.
We first recall that at the first-order level, the polygon $\mathcal{P}_n$ converges 
to a deterministic shape, namely
\begin{itemize}
\item
The distance $d(A_n,S_n)/n$ converges in distribution to 0 (and likewise for $d(B_n,E_n)/n$ etc.)
\item
The part of the boundary of $\mathcal{P}_n$ between $A_n$ and $B_n$, converges, after
renormalization, to an arc of parabola (and likewise  between $B_n$ and $C_n$ etc.)
\end{itemize}
A formal statement can be found in Barany \cite{Barany:limitshape}. 

Theorem 2 gives the second-order asymptotics.
We shall need an estimate on the number of convex lattice chains from $(0,0)$ to $(n+\delta_n,n+\delta'_n)$, contained in $[0, n+\delta_n]\times [0, n+\delta'_n]$.
The logarithm of this number is given by (see the proof of Lemma 2.1 in Bureaux-Enriquez \cite{Bureaux:polygons})
$$
  \frac{\zeta(3)}{\zeta(2)}\frac1{\beta_1\beta_2}+\frac{1}{2i\pi} \int_{1+\delta-i\infty}^{1+\delta+i\infty} \frac{\Gamma(s) \zeta(s+1)}{\zeta(s)}\chi(s;\bbeta)\,ds.
$$
where
$$
  \chi(s;\bbeta):=
  \sum_{\bv \in \Z_+^2\setminus\{0\}}\frac{1}{(\bbeta \cdot \bv)^s},
  \qquad \Re(s) > 2.
$$
and 
$\bbeta=(\beta_1, \beta_2)=(1/(n+\delta_n), 1/(n+\delta'_n)) $.
Moreover, by the same lemma, for all nonnegative integers $k_1, k_2$, all $\epsilon>0$ and all
$\bbeta=(\beta_1, \beta_2)\in (0,+\infty)^2$, such that
$\epsilon < \frac{\beta_1}{\beta_2}<\frac1{\epsilon}$,
$$
  \frac{\partial^{k_1+k_2}}{\partial\beta_1^{k_1}\partial\beta_2^{k_2}}\log Z(\beta_1, \beta_2)
  \underset{\bbeta\to0}\sim
  (-1)^{k_1+k_2}\frac{\zeta(3)}{\zeta(2)}\frac{k_1!k_2!}{\beta_1^{k_1+1}\beta_2^{k_2+1}}.
$$
Bureaux and Enriquez use these estimates to study the function $\chi$ near $(0,0)$ along the diagonal, that is, in  the case when $\delta_n=\delta'_n=0$. However, they point out that this can be extended to a more general setting. The case when $\delta_n$ and $\delta'_n$ are $o(n)$ corresponds  to studying $\chi$ in the neighbourhood of $(0,0)$ and near the diagonal. Using exactly the same arguments as in Lemma 2.2 leads to the following estimate. There exist $C, K>0$ such that
for all sequences $(\delta_n), (\delta'_n)$ satisfying $\delta_n=o(n)$,
$\delta'_n=o(n)$,
the number of convex lattice chains from $(0,0)$ to $(n+\delta_n,n+\delta'_n)$,
contained in $[0, n+\delta_n]\times [0, n+\delta'_n]$ and without vertical steps (it is easy to see that this last additional condition on vertical steps does not change the form of the estimate) is given by
\begin{equation}\label{chain}
\exp(C(n+\delta_n)^{1/3}(n+\delta'_n)^{1/3}+H(n+\delta_n,n+\delta'_n)+K\log n+o(1))
\end{equation}
where $H$ is a corrective term such that, uniformly over all $s, t\in \R$,
\begin{equation}\label{H}
|H(n+sn^{2/3},n+tn^{2/3})|=o(n^{1/3})
\end{equation}
On the other hand, by the same arguments, for fixed reals $s,t>0$,
the number of convex lattice chains from $(0,0)$ to $(sn,tn)$,
contained in $[0, sn]\times [0, tn]$ and such that
$(sn,tn)$ is the only point in the chain whose $x$-coordinate
is $sn$ is given by
\begin{equation}\label{chain2}
\exp(C(sn)^{1/3}(nt)^{1/3}+o(n^{1/3})
\end{equation}

We are now ready to prove Theorem 2.

(i) We use the same arguments as in \cite{Barany:limitshape}.
A convex lattice polygon can be seen as the union of 4 lattice convex chains.
In particular the number of polygons of this kind contained in
$[-n,n]^2$ is at least the number of such polygons
satisfying $A_n=S_n, B_n=E_n$ etc. which, according to \eqref{chain}, has the form
\begin{equation}\label{centre}
  \exp(4Cn^{2/3}+4K\log n+4H(n,n)+o(1))
  \end{equation}

If $c>0$, if we want to count the number of convex lattice polygons
contained in the rectangle $R_n:=[-n,n]\times[-(n-cn^{(1/3)+\delta}),n]$, we can choose the points
$A_n,B_n,C_n,D_n$ and count the lattice chains in-between. Comparing \eqref{chain} and
\eqref{chain2}, we see that this number is maximized when $A_n,B_n,C_n,D_n$
lie near the middle of the segments of $R_n$. For such a choice of $A_n,B_n,C_n,D_n$,
this  number of convex lattice polygons is bounded above by
$$\exp(4Cn^{1/3}(n-(c/2)n^{2/3})^{1/3}+4H(n, n-(c/2)n^{2/3})+4K\log n+o(1))$$
On the other hand, the number of choices of  $(A_n,B_n,C_n,D_n)$ is bounded $(2n)^4$, so that
the number of convex lattice polygons
contained in $R_n$ is bounded by
$$
\exp(4Cn^{2/3}-\frac{2}{3}Ccn^{1/3} + 4H(n, n-(c/2)n^{2/3})+(4K+16)\log n+o(1))$$
Using \eqref{H}, we obtain
\begin{equation}\label{need}
\frac{4H(n,n)}{-\frac{2}{3}Ccn^{1/3} + 4H(n, n-(c/2)n^{2/3})}\to \infty
 \end{equation}
as $n\to\infty$, which entails that that the probability that a
random, uniform  convex lattice polygon
contained $[-n,n]^2$ lies in $R_n$ tends to 0 as $n\to\infty$.
So for every $\varepsilon$, with probability going to 1, $n^{-2/3}|Y(A_n- S_n)|<\varepsilon$
and we have the same estimates for $B_n, C_n, D_n$.

(ii) From (i), we know that there exists a sequence $(D_n)$ of integers such that $D_n=o(n^{2/3})$
and that with probability going to 1, $|Y(A_n)+n|<D_n$, $|X(B_n)-n|<D_n$ etc.
Let now $(\delta_n,\delta'_n,\delta''_n,\delta'''_n)$ be sequences of integers such that for each $n$,
$0<\delta_n<D_n$ etc. We want to study the law of 
$$n^{-2/3}(X(A_n- S_n), Y(B_n- E_n),-X(C_n- N_n)-Y(D_n- W_n))$$
conditionally on the event
$${\bf E}_n=\{Y(A_n- S_n)=\delta_n, X(B_n- E_n)=\delta'_n,Y(C_n- N_n)=\delta''_n,
X(D_n- W'_n)=\delta'''_n\}$$
Fix $r,s,t,u \in\R$. For simplicity, we shall omit the integer part notation in the sequel.
We want to estimate the number of 4-tuple of convex chains such that
\begin{eqnarray*}
&&X(A_n- S_n)=rn^{2/3}, Y(A_n- S_n)=\delta_n\\
&&Y(B_n- E_n)=sn^{2/3}, X(B_n- E_n)=\delta'_n
\end{eqnarray*}
etc.
Let $L_n(r,s,t,u)$ denote the logarithm of the number of such 4-tuples of
chains. The number of chains going from $A_n$ to $B_n$ is the same as the number of chains
from $(0,0)$ to $(n-\delta'_n-rn^{2/3},n-\delta_n+sn^{2/3})$. According to \eqref{chain}, this leads to
\begin{eqnarray*}
  \L_n(r,s,t,u)&=&Cn^{2/3}\left(\left(1-\frac{r}{3n^{1/3}}-\frac{2r^2}{9n^{2/3}}-\frac{\delta'_n}{n}
                   +O(1/n)\right)
    \left(1+\frac{s}{3n^{1/3}}-\frac{2s^2}{9n^{2/3}}-\frac{\delta_n}{n}+O(1/n))\right)\right)\\
&&+H\left(n-rn^{2/3}-\delta'_n,n+sn^{2/3}-\delta_n \right)+k\log n +\ldots
\end{eqnarray*}
where we only wrote the term corresponding to the chain from $A_n$ to $B_n$ but of course, there are
3 other terms. Summing up, we see that terms of the form $r/n^{1/3}$ cancel out and we get
\begin{eqnarray*}
  L_n(r,s,t,u)&=&Cn^{2/3}\left(4+ 2\left(\frac{\delta_n}{n}+\frac{\delta'_n}{n}+\frac{\delta''_n}{n}
                          +\frac{\delta'''_n}{n}\right)
                      -\frac{4}{9}\left(\frac{r^2+s^2+t^2+u^2}{n^{2/3}}\right)
                         -\frac{1}{9}\left(\frac{rs+st+tu+ur}{n^{2/3}}\right)\right)\\
                      &&+4K \log n+H\left(n-rn^{2/3}-\delta'_n,n+sn^{2/3}-\delta_n \right)+\ldots
                         +O(1/n)
\end{eqnarray*}
where again, in the last line, we have 3 additional terms involving the function $H$.
In particular, we get
$$
\frac{L_n(r,s,t,u)}{Cn^{2/3}}=\frac{L_n(0,0,0,0)}{Cn^{2/3}}
-\frac{1}{18}[(r-s)^2+(s-t)^2+(t-u)^2+(u-r)^2]-\frac{1}{3}
[r^2+s^2+t^2+u^2]+o(1)
$$
using the following fact that follows from \eqref{H}:
\begin{equation}\label{need1}
|H\left(n-rn^{2/3}-\delta'_n,n+sn^{2/3}-\delta_n \right)-
H\left(n-\delta'_n,n-\delta_n \right)|=o(n^{2/3})
\end{equation}
etc.
This is true for all sequences $(\delta_n,\delta'_n,\delta''_n,\delta'''_n)$, which completes the
proof.

(iii) This part is a reformulation of Proposition 6 in Section \ref{section_donsker}. In Theorem 2 (iii), we are looking at the southeast arc of the polygon whereas Proposition 6 deals with what would be the northwest arc. Hence  the identification
$$\mu_{r,s}(t)=\begin{pmatrix}
0 & 1\\
1 & 0
\end{pmatrix}\nu_{-r,s}(t)$$ 
The matrix $Q(t)$ in Theorem 2 (iii) corresponds to $P(t)D(t)$ in Proposition 6. Finally, the matrix $\mathcal{O}_t$ mentioned in the comments after the statement of Theorem 2 corresponds to the orthogonal matrix appearing in \eqref{orthogonal}.

%% file: Parties/Development/8acknowledgment.tex
We would like to thank Andrea Sportiello for useful discussions on the topic.